\documentclass[a4paper]{article}
\usepackage{graphicx}
\usepackage{amsmath}
\usepackage{amsfonts}
\usepackage{amssymb}
\newtheorem{theorem}{Theorem}[subsection]

\newtheorem{claim}[theorem]{Rule}

\newtheorem{definition}[theorem]{Definition}
\newtheorem{example}[theorem]{Example}

\newtheorem{lemma}[theorem]{Lemma}

\newtheorem{proposition}[theorem]{Proposition}
\newtheorem{remark}[theorem]{Remark}

\newenvironment{proof}[1][Proof]{\textbf{#1.} }{\ \rule{0.5em}{0.5em}}

\begin{document}

\title{An algorithm for computing the global basis of an irreducible $U_{q}(sp_{2n})$-module}
\author{C\'{e}dric Lecouvey\\lecouvey@math.unicaen.fr}
\date{}
\maketitle
\begin{abstract}
We describe a simple algorithm for computing the canonical basis of any
irreducible finite-dimensional $U_{q}(sp_{2n})$-module.
\end{abstract}

\section{Introduction}

The quantum algebra $U_{q}(\frak{g})$ associated to a semisimple Lie algebra
$\frak{g}$ is the $q$-analogue introduced by Drinfeld and Jimbo of its
universal enveloping algebra $U(\frak{g})$. Kashiwara \cite{Ka} and Lusztig
\cite{Lut} have discovered a distinguished basis of $U_{q}^{-}(\frak{g})$
which projects onto a global crystal basis (Kashiwara) or canonical basis
(Lusztig) of each simple finite-dimensional $U_{q}(\frak{g})$-module. When $q$
tends to $1$ this basis yields in particular a canonical basis of the
corresponding $U(\frak{g})$-module.

In this article, we restrict ourselves to the case $g=sp_{2n}$. Denote by
$\{\Lambda_{1},...,\Lambda_{n}\}$ the set of fundamental weights, by $P$ the
weight lattice and by $P_{+}$ the set of dominant weights of $sp_{2n}.$ Then
for each $\lambda=\sum\lambda_{i}\Lambda_{i}\in P_{+}$ there exists a unique
irreducible finite-dimensional $U_{q}(sp_{2n})$-module $V(\lambda)$ of highest
weight $\lambda$. The aim of this article is to describe a simple algorithm
for computing the global crystal basis of $V(\lambda)$. This algorithm will be
a generalization of the one described by Leclerc and Toffin in \cite{L-T} for
the irreducible $U_{q}(sl_{n})$-modules. In the case of $U_{q}(sp_{2n})$, an
algorithm was only known for the fundamental modules $V(\Lambda_{p})$,
$p=1,...,n$ \cite{Ma}. We note that Zelevinsky and Retakh \cite{ZR} have
described for every irreducible finite-dimensional $U(sp_{4})$-module a
so-called good basis.\ This basis is the specialization at $q=1$ of the dual
of the canonical basis.

Our method is as follows. First we realize $V(\Lambda_{p})$ as a
subrepresentation of a $U_{q}(sp_{2n})$-module $W(\Lambda_{p})$ whose basis
$\{v_{C}\}$ has a natural indexation in terms of column shaped Young tableaux.
This representation $W(\Lambda_{p})$ may be regarded as a $q$-analogue of the
$p$-th exterior power of the vector representation of $sp_{2n}$. Then we give
explicit formulas for the expansion of the global crystal basis of
$V(\Lambda_{p})$ on the basis $\{v_{C}\}$. Next we embed $V(\lambda)$ in the
tensor product%
\[
W(\lambda)=W(\Lambda_{1})^{\bigotimes\lambda_{1}}%
{\textstyle\bigotimes}
\cdot\cdot\cdot%
{\textstyle\bigotimes}
W(\Lambda_{n})^{\bigotimes\lambda_{n}}%
\]
The tensor product of the crystal bases of the $W(\Lambda_{p})$'s is a natural
basis $\{v_{\tau}\}$ of $W(\lambda)$ indexed by combinatorial objects $\tau$
called tabloids. Then we obtain an intermediate basis of $V(\lambda)$ fixed by
the involution $q\longmapsto q^{-1}$ and such that the transition matrix from
this basis to the global crystal basis of $V(\lambda)$ is unitriangular.
Finally we compute the expansion of the canonical basis on the basis
$\{v_{\tau}\}$ via an elementary algorithm. We give as an example the matrix
associated to the expansion of the global basis of a weight space of the
$U_{q}(sp_{6})$-module $V(\lambda)$ with $\lambda=\Lambda_{1}+\Lambda
_{2}+2\Lambda_{3}.$

\section{Background}

In this section we briefly review the basic facts that we shall need
concerning the representation theory of $U_{q}(sp_{2n})$ and the notions of
crystal basis and canonical basis of a $U_{q}(sp_{2n})$-module. The reader is
referred to \cite{Ka2}, \cite{Ka1}, \cite{Ch-Pr} and \cite{Jan} for more details.

\subsection{The quantum enveloping algebra $U_{q}(sp_{2n})$}

Recall the Dynkin diagram of $sp_{2n}$:%
\[
\overset{1}{\circ}-\overset{2}{\circ}-\overset{3}{\circ}\cdot\cdot
\cdot\overset{n-2}{\circ}-\overset{n-1}{\circ}\Longleftarrow\overset{n}{\circ
}.
\]
Accordingly, the Cartan matrix $A=(a_{i,j})$ of $sp_{2n}$ is:%
\[
\left(
\begin{array}
[c]{cccccc}%
\text{ \ }2 & -1 &  &  &  & \\
-1 & \text{ \ }2 & -1 &  &  & \\
& -1 & . & . &  & \\
&  & . & . & . & \\
&  &  & . & \text{ \ }2 & -1\\
&  &  &  & -2 & \text{ \ }2
\end{array}
\right)
\]
with rows and columns indexed by $\{1,...,n\}$. Given a fixed indeterminate
$q$ set%
\begin{gather*}
q_{i}=\left\{
\begin{tabular}
[c]{l}%
$q$ if $i\neq n$\\
$q^{2}$ if $i=n$%
\end{tabular}
\right.  \text{,}\\
\lbrack n]_{i}=\frac{q_{i}^{n}-q_{i}^{-n}}{q_{i}-q_{i}^{-1}}\text{ and
}[n]_{i}!=[n]_{i}[n-1]_{i}\cdot\cdot\cdot\lbrack1]_{i}.
\end{gather*}
The quantized enveloping algebra $U_{q}(sp_{2n})$ is the associative algebra
over $\mathbb{C}(q)$ generated by $e_{i},f_{i},q^{h_{i}},q^{-h_{i}}$
$i=1,...,n$, subject to the relations:%
\begin{gather}
q^{h_{i}}q^{-h_{i}}=q^{-h_{i}}q^{h_{i}}=1,\label{D1}\\
q^{\pm h_{i}}q^{\pm h_{j}}=q^{\pm h_{j}}q^{\pm h_{i}},\label{D2}\\
q^{h_{i}}e_{j}q^{-h_{i}}=q^{a_{i,j}}e_{j},\label{D3}\\
q^{h_{i}}f_{j}q^{-h_{i}}=q^{-a_{i,j}}f_{j},\label{D4}\\
\lbrack e_{i},f_{i}]=\frac{t_{i}-t_{i}^{-1}}{q_{i}-q_{i}^{-1}}\delta
_{i,j}\text{ where }t_{i}=\left\{
\begin{tabular}
[c]{l}%
$q^{h_{i}}$ if $i\neq n$\\
$q^{2h_{n}}$ if $i=n$%
\end{tabular}
\right.  ,\label{D5}\\
\text{if }i\neq j\text{ }\overset{1-a_{j,i}}{\underset{k=0}{\sum}}%
(-1)^{k}e_{i}^{(k)}e_{j}e_{i}^{(1-a_{j,i}-k)}=\overset{1-a_{j,i}}%
{\underset{k=0}{\sum}}(-1)^{k}f_{i}^{(k)}f_{j}f_{i}^{(1-a_{j,i}-k)}=0
\label{D6}%
\end{gather}
where $e_{i}^{(m)}=e_{i}^{m}/[m]_{i}!$ and $f_{i}^{(m)}=f_{i}^{m}/[m]_{i}!$.

The subalgebra of $U_{q}(sp_{2n})$ generated by $e_{i},f_{i},q^{h_{i}%
},q^{-h_{i}}$ $i=1,...,n-1$ is isomorphic to $U_{q}(sl_{n})$, the quantum
enveloping algebra of $sl_{n}$.

The representation theory of $U_{q}(sp_{2n})$ is closely parallel to that of
its classical counterpart $U(sp_{2n})$. The weight lattice $P$ of
$U_{q}(sp_{2n})$ is the $\mathbb{Z}$-lattice generated by the fundamental
weights $\Lambda_{1},...,\Lambda_{n}.$ We denote by $P_{+}$ the set of
dominant weights of $U_{q}(sp_{2n})$ i.e. those of the form $\lambda
_{1}\Lambda_{1}+\cdot\cdot\cdot+\lambda_{n}\Lambda_{n}$ where $\lambda
_{1},...,\lambda_{n}\in\mathbb{N}$. Let $M$ be a $U_{q}(sp_{2n})$-module. For
every $\mu\in P$\ the subspace%
\[
M_{\mu}=\{v\in M,\text{ }q^{h_{i}}v=q^{<h_{i},\mu>}v,\text{ }i=1,...,n\}
\]
is the weight space of weight $\mu$ of $M$. A vector $v_{0}$ is said to be of
highest weight when $e_{i}(v_{0})=0$ for $i=1,...,n$. If $M$ is a
finite-dimensional irreducible module, $M$ is a highest weight module that is,
contains a highest weight vector $v_{\lambda}$ of weight $\lambda$ such that
$M=U_{q}(sp_{2n})v_{\lambda}$. Then $\dim M_{\lambda}=1$ and $\lambda$ is a
dominant weight. Conversely, for each dominant weight $\lambda\in P_{+},$
there is a unique finite-dimensional module with highest weight $\lambda$.\ We
denote it by $V(\lambda)$ and we write $v_{\lambda}$ for a fixed highest
weight vector.

Given two $U_{q}(sp_{2n})$-modules $M$ and $N$, we can define a structure of
$U_{q}(sp_{2n})$-module on $M%
{\textstyle\bigotimes}
N$ by putting:%
\begin{gather}
q^{h_{i}}(u%
{\textstyle\bigotimes}
v)=q^{h_{i}}u%
{\textstyle\bigotimes}
q^{h_{i}}v,\label{tensor1}\\
e_{i}(u%
{\textstyle\bigotimes}
v)=e_{i}u%
{\textstyle\bigotimes}
t_{i}^{-1}v+u%
{\textstyle\bigotimes}
e_{i}v,\label{tensor2}\\
f_{i}(u%
{\textstyle\bigotimes}
v)=f_{i}u%
{\textstyle\bigotimes}
v+t_{i}u%
{\textstyle\bigotimes}
f_{i}v. \label{tensor3}%
\end{gather}

\subsection{Crystal basis for $U_{q}(sp_{2n})$-modules \label{subsec_crystal}}

Let $M$ be a finite-dimensional $U_{q}(sp_{2n})$-module. Fix $i\in
\{1,...,n\}$. Let $u\in M$ and suppose that $u\in M_{\mu}$. Then $u$ can be
written uniquely as a finite sum $\underset{k}{\sum}f_{i}^{(k)}u_{k}$ where
$u_{k}\in M_{\mu+k\alpha_{i}}$ and $e_{i}u_{k}=0$. Kashiwara's operators
$\widetilde{e}_{i}$ and $\widetilde{f}_{i}$ are defined\ by:%
\begin{equation}
\widetilde{e}_{i}u=\sum f_{i}^{(k-1)}u_{k}\text{,\ \ \ \ \ \ \ \ \ }%
\widetilde{f}_{i}u=\sum f_{i}^{(k+1)}u_{k}\text{.} \label{action_tilda}%
\end{equation}
Denote by $A$ the subalgebra of $\mathbb{C}(q)$ consisting of the rational
functions without pole at $q=0$. Let $L$ be a free $A$-submodule of $M$ such
that $L%
{\textstyle\bigotimes}
\mathbb{C(}q)=M$ and $B$ a basis of the $\mathbb{Q}$-vector space $L/qL$.
Write $\pi$ for the canonical projection
\[
L\overset{\pi}{\rightarrow}L/qL
\]
Set $L_{\mu}=L\cap M_{\mu}$ and $B_{\mu}=B\cap L_{\mu}/qL_{\mu}$. Then $(L,B)$
is a crystal basis of $M$ at $q=0$ if the following conditions hold:%

\begin{equation}
\left\{
\begin{tabular}
[c]{l}%
$(\mathrm{i})$ $L=\underset{\mu\in P}{%
{\textstyle\bigoplus}
}L_{\mu}$ and $B=\underset{\mu\in P}{\cup}B_{\mu},$\\
$(\mathrm{ii})$ $\widetilde{e}_{i}L\subset L\text{ and }\widetilde{f}%
_{i}L\subset L,$\\
$(\mathrm{iii})$ $\widetilde{e}_{i}B\subset B\cup\{0\}\text{ and }%
\widetilde{f}_{i}B\subset B\cup\{0\},$\\
$(\mathrm{iv})\text{ for }b_{1},b_{2}\in B\text{ and }i\in\{1,...,n\},\text{
}\widetilde{e}_{i}b_{1}=b_{2}\Longleftrightarrow b_{1}=\widetilde{f}_{i}%
b_{2}.$%
\end{tabular}
\right.  \label{def_C_B}%
\end{equation}
Note that the action of $\widetilde{e}_{i}$ and $\widetilde{f}_{i}$ on $L/qL$
is well defined because of $(\mathrm{ii})$.\ Kashiwara \cite{Ka0} has proved
that every finite-dimensional $U_{q}(sp_{2n})$-module $M$ has a crystal basis.
Moreover if $M=V(\lambda)$ is simple, this basis is unique up to an overall
scalar factor. We shall denote it by $(L(\lambda),B(\lambda))$.

If $(L,B)$ and $(L^{\prime},B^{\prime})$ are crystal bases of the
finite-dimensional $U_{q}(sp_{2n})$-modules $M$ and $M^{\prime}$, then $(L%
{\textstyle\bigotimes}
L^{\prime},$ $B%
{\textstyle\bigotimes}
B^{\prime})$ with $B\bigotimes B^{\prime}=\{b\bigotimes b^{\prime};$ $b\in
B,b^{\prime}\in B^{\prime}\}$ is a crystal basis of $M%
{\textstyle\bigotimes}
M^{\prime}$. The action of $\widetilde{e}_{i}$ and $\widetilde{f}_{i}$ on
$B\bigotimes B^{\prime}$ is given by:%

\begin{align}
\widetilde{f_{i}}(u%
{\textstyle\bigotimes}
v)  &  =\left\{
\begin{tabular}
[c]{c}%
$\widetilde{f}_{i}(u)\bigotimes v$ if $\varphi_{i}(u)>\varepsilon_{i}(v)$\\
$u\bigotimes\widetilde{f}_{i}(v)$ if $\varphi_{i}(u)\leq\varepsilon_{i}(v)$%
\end{tabular}
\right. \label{TENS1}\\
&  \text{and}\nonumber\\
\widetilde{e_{i}}(u%
{\textstyle\bigotimes}
v)  &  =\left\{
\begin{tabular}
[c]{c}%
$u\bigotimes\widetilde{e_{i}}(v)$ if $\varphi_{i}(u)<\varepsilon_{i}(v)$\\
$\widetilde{e_{i}}(u)\bigotimes v$ if$\varphi_{i}(u)\geq\varepsilon_{i}(v)$%
\end{tabular}
\right.  \label{TENS2}%
\end{align}
where $\varepsilon_{i}(u)=\max\{k;\widetilde{e}_{i}^{k}(u)\neq0\}$ and
$\varphi_{i}(u)=\max\{k;\widetilde{f}_{i}^{k}(u)\neq0\}$.

The set $B$ may be endowed with a combinatorial structure called the crystal
graph of $M$. Crystal graphs are oriented colored graphs with colors
$i\in\{1,...,n\}.$ An arrow $a\overset{i}{\rightarrow}b$ means that
$\widetilde{f}_{i}(a)=b$ and $\widetilde{e}_{i}(b)=a$. The decomposition of
$M$ into its irreducible components is reflected into the decomposition of $B$
into its connected components. The crystal graphs of two isomorphic
irreducible components are isomorphic as oriented colored graphs. A vertex
$v^{0}\in B$ satisfying $\widetilde{e}_{i}(v^{0})=0$ for $i\in\{1,..,n\}$ is
called a highest weight vertex. The crystal $B(\lambda)$ contains a unique
highest weight vertex.

The end of this section is devoted to Kashiwara-Nakashima's combinatorial
description of the crystal graphs of the finite-dimensional irreducible
$U_{q}(sp_{2n})$-modules \cite{KN}. It is based on the notion of symplectic
tableaux analogous to Young tableaux for type $A$. In the sequel we use De
Concini's version of these tableaux which is equivalent to
Kashiwara-Nakashima's one \cite{Sh}.

Let us consider the totally ordered alphabet
\[
\mathcal{C}_{n}=\{1<\cdot\cdot\cdot<n<\overline{n}<\cdot\cdot\cdot
<\overline{1}\}.
\]
For each letter $x\in\mathcal{C}_{n}$ we denote by pred$(x)$ the largest
letter $y$ such that $y<x$. A column on $\mathcal{C}_{n}$ is a Young diagram
$C$ of column shape filled from top to bottom by increasing letters of
$\mathcal{C}_{n}$. The height $h(C)$ of a column $C$ is the number of its
letters. Set $\mathbf{C}(n,h)$ for the set of columns of height $h$ on
$\mathcal{C}_{n}$ i.e. with letters in $\mathcal{C}_{n}$. The reading of the
column $C\in\mathbf{C}(n,h)$ is the word $\mathrm{w}(C)$ of $\mathcal{C}%
_{n}^{\ast}$ obtained by reading the letters of $C$ from top to bottom.\ We
will say that a column $C$ contains the pair $(z,\overline{z})$ when $C$
contains the unbarred letter $z\leq n$ and the barred letter $\overline{z}%
\geq\overline{n}$.\ Let $C_{1}$ and $C_{2}$ be two columns. We will write
$C_{1}\leq C_{2}$ when $h(C_{1})\geq h(C_{2})$ and the rows of the tableau
$C_{1}C_{2}$ weakly increase.

\begin{definition}
Let $C$ be a column and $I_{C}=\{z_{1}>\cdot\cdot\cdot>z_{r}\}$ the set of
unbarred letters $z$ such that the pair $(z,\overline{z})$ occurs in $C$. The
column $C$ is admissible when there exists a set of unbarred letters
$J_{C}=\{t_{1}>\cdot\cdot\cdot>t_{r}\}\subset\mathcal{C}_{n}$ such that:

\begin{itemize}
\item $t_{1}$ is the greatest letter of $\mathcal{C}_{n}$ satisfying:
$t_{1}<z_{1},t_{1}\notin C$ and $\overline{t_{1}}\notin C,$

\item  for $i=2,...,r$, $t_{i}$ is the greatest letter of $\mathcal{C}_{n}$
satisfying: $t_{i}<\min(t_{i-1,}z_{i}),$ $t_{i}\notin C$ and $\overline{t_{i}%
}\notin C.$
\end{itemize}

In this case we write:

\begin{itemize}
\item $rC$ for the column obtained from $C$ by changing $\overline{z}_{i}$
into $\overline{t}_{i}$ for each letter $z_{i}\in I_{C},$

\item $lC$ for the column obtained from $C$ by changing $z_{i}$ into $t_{i}$
for each letter $z_{i}\in I_{C}.$
\end{itemize}
\end{definition}

A column $C$ on $\mathcal{C}_{n}$ may be non admissible. For $C=%
\begin{tabular}
[c]{|l|}\hline
$\mathtt{2}$\\\hline
$\mathtt{3}$\\\hline
$\mathtt{\bar{3}}$\\\hline
$\mathtt{\bar{1}}$\\\hline
\end{tabular}
$ it is impossible to find a letter $t<3$ such that $t\notin C$ and
$\overline{t}\notin C.\;$We write $\mathbf{Ca}(n,h)$ for the set of admissible
columns of height $h$ on $\mathcal{C}_{n}$.

Notice that the condition $t_{i}<\min(t_{i-1},z_{i})$ for $i=2,...,r$ of the
above definition can be replaced by the condition $t_{i}<z_{i}$ and
$t_{i}\notin\{t_{1},...,t_{i-1}\}$ for $i=2,...,r$.\ Indeed, $t_{i}>t_{i-1}$
contradicts the fact that $t_{i-1}$ is maximal. Then $J_{C}$ may be regarded
as the set of $r$ maximal unbarred letters $\{t_{1},...,t_{r}\}$ such that
$t_{i}<z_{i}$ and $\{t_{i},\overline{t}_{i}\}\cap C=\emptyset$.\ 

Similarly to the type $A$ case, we can associate to each dominant weight
$\lambda=\overset{n}{\underset{i=1}{\sum}}\lambda_{i}\Lambda_{i}$ a Young
diagram $Y(\lambda)$ having $\lambda_{i}$ columns of height $i,$ $i=1,...,n$.
By definition, a symplectic tableau $T$ of shape $\lambda$ is a filling of
$Y(\lambda)$ by letters of $\mathcal{C}_{n}$ satisfying the following conditions:

\begin{itemize}
\item  the columns $C_{i}$ of $T=C_{1}\cdot\cdot\cdot C_{s}$ are admissible,

\item  for $i=1,...,s-1:rC_{i}\leq lC_{i+1}.$
\end{itemize}

The set of symplectic tableaux of shape $\lambda$ will be denoted
$\mathbf{ST}(n,\lambda)$. If $T=C_{1}C_{2}\cdot\cdot\cdot C_{r}\in
\mathbf{ST}(n,\lambda)$, the reading of $T$ is the word $\mathrm{w}%
(T)=\mathrm{w}(C_{r})\cdot\cdot\cdot\mathrm{w}(C_{2})\mathrm{w}(C_{1})$.\ From
\cite{KN} and \cite{Sh} we deduce the

\begin{theorem}
\ \ \ \ \ \label{TH_KN}

\textrm{(i)}: The vertices of $B(\Lambda_{p})$ are in one-to-one
correspondence with the readings of admissible columns of height $p.$

\textrm{(ii)}: The vertices of $B(\lambda)$ are in one-to-one correspondence
with the readings of the symplectic tableaux of shape $\lambda$.
\end{theorem}

More precisely Kashiwara and Nakashima realize $V(\lambda)$ into a tensor
power $V(\Lambda_{1})^{\bigotimes l}$ of the vector representation whose
crystal graph is:%
\begin{equation}
1\overset{1}{\rightarrow}2\cdot\cdot\cdot\cdot\rightarrow n-1\overset
{n-1}{\rightarrow}n\overset{n}{\rightarrow}\overline{n}\overset{n-1}%
{\rightarrow}\overline{n-1}\overset{n-2}{\rightarrow}\cdot\cdot\cdot
\cdot\rightarrow\overline{2}\overset{1}{\rightarrow}\overline{1}.
\label{crystal_V}%
\end{equation}
Let us identify the vertices of the crystal graph $G_{n}=\underset{l}{%
{\textstyle\bigoplus}
}B(\Lambda_{1})^{\bigotimes l}$ with the words on $\mathcal{C}_{n}$.$\;$The
weight of the vertex $b\in G_{n}$ is defined by%
\begin{equation}
\mathrm{wt}(b)=\underset{i=1}{\overset{n}{\sum}}(\varphi_{i}(b)-\varepsilon
_{i}(b))\Lambda_{i}. \label{def_wt}%
\end{equation}
$B(\Lambda_{p})$ can then be identified with the connected component of
$G_{n}$ whose vertex of highest weight is the reading of the column
\[
C_{p}^{0}=%
\begin{tabular}
[c]{|l|}\hline
$1$\\\hline
$2$\\\hline
$\cdot$\\\hline
$\cdot$\\\hline
$p$\\\hline
\end{tabular}
\text{ .}%
\]
In this identification, the vertices of $B(\Lambda_{p})$ are the readings of
the admissible columns of height $p$.\ If $\lambda=\underset{p=1}{\overset
{n}{\sum}}\lambda_{p}\Lambda_{p}$, $B(\lambda)$ is identified with the
connected component whose highest weight vertex is the reading of the
symplectic tableau $T_{\lambda}$ containing $\lambda_{p}$ columns $C_{p}^{0}$
for $p=1,...,n$. Then the vertices of $B(\lambda)$ are the readings of the
symplectic tableaux of shape $\lambda$.

Using Formulas (\ref{TENS1}) and (\ref{TENS2}) we obtain a simple rule to
compute the action of $\widetilde{e}_{i}$ and $\widetilde{f}_{i}$ on $w\in
G_{n}$ that we will use in Section \ref{sec_G_lambda}. Consider the subword
$w_{i}$ of $w$ containing only the letters $\overline{i+1},\overline{i}%
,i,i+1$. Then encode in $w_{i}$ each letter $\overline{i+1}$ or $i$ by the
symbol $+$ and each letter $\overline{i}$ or $i+1$ by the symbol $-$. Because
$\widetilde{e}_{i}(+-)=\widetilde{f}_{i}(+-)=0$ in $B(\Lambda_{1})%
{\textstyle\bigotimes}
B(\Lambda_{1})$ the factors of type $+-$ may be ignored in $w_{i}.$ So we
obtain a subword $w_{i}^{(1)}$ in which we can ignore all the factors $+-$ to
construct a new subword $w_{i}^{(2)}$ etc... Finally we obtain a subword
$\rho(w)$ of $w$ of type%
\[
\rho(w)=-^{r}+^{s}.
\]
Then we have the

\begin{claim}
\ \ \ \ \ \ \label{+-}

\begin{itemize}
\item  If $r>0,$ $\widetilde{e}_{i}(w)$ is obtained by changing the rightmost
symbol $-$ of $\rho(w)$ into its corresponding symbol $+$ (i.e. $i+1$ into $i$
and $\overline{i}$ into $\overline{i+1}$) the others letters of $w$ being
unchanged. If $r=0,$ $\widetilde{e}_{i}(w)=0$

\item  If $s>0,$ $\widetilde{f}_{i}(w)$ is obtained by changing the leftmost
symbol $+$\ of $\rho(w)$ into its corresponding symbol $-$ (i.e. $i$ into
$i+1$ and $\overline{i+1}$ into $\overline{i}$) the others letters of $w$
being unchanged. If $s=0,\widetilde{f}_{i}(w)=0.$
\end{itemize}
\end{claim}

\subsection{Canonical bases for $U_{q}(sp_{2n})$-modules}

In the sequel we identify $B(\Lambda_{p})$ to $\{\mathrm{w}(C);$
$C\in\mathbf{Ca}(n,p)\}$ and $B(\lambda)$ to $\{\mathrm{w}(T);$ $T\in
\mathbf{ST}(n,\lambda)\}$. Denote by $F\mapsto\overline{F}$ the involution of
$U_{q}(sp_{2n})$ defined as the ring automorphism satisfying%
\[
\overline{q}=q^{-1},\text{ \ \ }\overline{q^{h_{i}}}=q^{-h_{i}},\text{
\ \ }\overline{e_{i}}=e_{i},\text{ \ \ }\overline{f_{i}}=f_{i}\text{ \ \ \ for
}i=1,...,n.
\]
Writing each vector $v$ of $V(\lambda)$ in the form $v=Fv_{\lambda}$ where
$F\in U_{q}(sp_{2n})$, we obtain an involution of $V(\lambda)$ defined by%
\[
\overline{v}=\overline{F}v_{\lambda}.
\]
Let $U_{\mathbb{Q}}^{-}$ be the subalgebra of $U_{q}(sp_{2n})$ generated over
$\mathbb{Q}[q,q^{-1}]$ by the $f_{i}^{(k)}$ and set $V_{\mathbb{Q}}%
(\lambda)=U_{\mathbb{Q}}^{-}v_{\lambda}$. We can now state:

\begin{theorem}
\label{TH_K_2}(Kashiwara)

There exists a unique $\mathbb{Q[}q,q^{-1}]$-basis $\{G(T);$ $T\in
\mathbf{ST}(n,\lambda)\}$ of $V_{\mathbb{Q}}(\lambda)$ such that:%
\begin{gather}
G(T)\equiv\mathrm{w(}T)\text{ }\mathrm{mod}\text{ }qL(\lambda
),\label{cond_cong}\\
\overline{G(T)}=G(T). \label{cond_invo}%
\end{gather}
\end{theorem}

\noindent This basis is called the lower global (or canonical) basis of
$V(\lambda)$, and our aim is to calculate it.

\section{Fundamental modules}

\subsection{Marsh's Algorithm\label{subsec_marsh}}

We review Marsh's algorithm for computing the global basis of $V(\Lambda_{p})$
\cite{Ma}. Let $\mathrm{w(}C)\in B(\Lambda_{p}).$ A letter $x\in C$ is said to
be movable if pred$(x)\notin C$. We define a path in $B(\Lambda_{p})$ joining
$\mathrm{w(}C)$ to $\mathrm{w(}C_{p}^{0})$. If $C\neq C_{p}^{0},$ let $z$ be
the lowest movable letter of $C.$ Then we compute a new column $C_{1}$ as
follows:%
\begin{align*}
\text{\textrm{(i)}}  &  :\text{ if }z=i+1\text{, }\overline{i}\notin C\text{
or }\overline{i+1}\in C,\text{ }C_{1}=C-\{i+1\}+\{i\},\\
\text{\textrm{(ii)}}  &  :\text{ if }z=i+1\text{, }\overline{i}\in C\text{ and
}\overline{i+1}\notin C,\text{ }C_{1}=C-\{\overline{i},i+1\}+\{\overline
{i+1},i\},\\
\text{\textrm{(iii)}}  &  :\text{ if }z=\overline{i}>\overline{n}\text{,
}C_{1}=C-\{\overline{i}\}+\{\overline{i+1}\},\\
\text{\textrm{(iv)}}  &  :\text{ if }z=\overline{n}\text{, }C_{1}%
=C-\{\overline{n}\}+\{n\}.
\end{align*}

\begin{remark}
\label{rem_util}In case \textrm{(iii) }the letters $i$ and $\overline{i}$ can
not appear simultaneously in $C$.\ Otherwise by definition of $z$, the letters
$i,i-1,...,1$ would appear in $C$ and $C$ would be not admissible.
\end{remark}

We will have $\mathrm{w(}C_{1})=\widetilde{e}_{i}\mathrm{w(}C)$ in cases
\textrm{(i)} and \textrm{(iii)}, $\mathrm{w(}C_{1})=\widetilde{e}_{i}%
^{2}\mathrm{w(}C)$ in case \textrm{(ii)} and $\mathrm{w(}C_{1})=\widetilde
{e}_{n}\mathrm{w(}C)$ in case \textrm{(iv)}. So $C_{1}$ is an admissible
column. Write $\mathrm{w(}C_{1})=\widetilde{e}_{i_{1}}^{p_{1}}\mathrm{w(}C)$.
Next we compute similarly $C_{2}$ from $C_{1}$ and write $\mathrm{w(}%
C_{2})=\widetilde{e}_{i_{2}}^{p_{2}}\mathrm{w(}C_{1})$. Finally , after a
finite number of steps, we will reach $C_{p}^{0}$ and we will get
$\mathrm{w(}C_{p}^{0})=\widetilde{e}_{i_{r}}^{p_{r}}\cdot\cdot\cdot
\widetilde{e}_{i_{1}}^{p_{1}}\mathrm{w(}C)$, hence $\mathrm{w(}C)=\widetilde
{f}_{i_{1}}^{p_{1}}\cdot\cdot\cdot\widetilde{f}_{i_{r}}^{p_{r}}\mathrm{w(}%
C_{p}^{0})$.

\begin{example}
Let $C=%
\begin{tabular}
[c]{|l|}\hline
$\mathtt{\bar{2}}$\\\hline
$\mathtt{\bar{1}}$\\\hline
\end{tabular}
$ and $n=3$. We obtain:\vspace{0.5cm}

$\mathrm{w}%
\begin{tabular}
[c]{|l|}\hline
$\mathtt{\bar{2}}$\\\hline
$\mathtt{\bar{1}}$\\\hline
\end{tabular}
\overset{\widetilde{e}_{2}}{\rightarrow}\mathrm{w}%
\begin{tabular}
[c]{|l|}\hline
$\mathtt{\bar{3}}$\\\hline
$\mathtt{\bar{1}}$\\\hline
\end{tabular}
\overset{\widetilde{e}_{3}}{\rightarrow}\mathrm{w}%
\begin{tabular}
[c]{|l|}\hline
$\mathtt{3}$\\\hline
$\mathtt{\bar{1}}$\\\hline
\end{tabular}
\overset{\widetilde{e}_{2}}{\rightarrow}\mathrm{w}%
\begin{tabular}
[c]{|l|}\hline
$\mathtt{2}$\\\hline
$\mathtt{\bar{1}}$\\\hline
\end{tabular}
\overset{\widetilde{e}_{1}^{2}}{\rightarrow}\mathrm{w}%
\begin{tabular}
[c]{|l|}\hline
$\mathtt{1}$\\\hline
$\mathtt{\bar{2}}$\\\hline
\end{tabular}
\overset{\widetilde{e}_{2}}{\rightarrow}\mathrm{w}%
\begin{tabular}
[c]{|l|}\hline
$\mathtt{1}$\\\hline
$\mathtt{\bar{3}}$\\\hline
\end{tabular}
\overset{\widetilde{e}_{3}}{\rightarrow}\mathrm{w}%
\begin{tabular}
[c]{|l|}\hline
$\mathtt{1}$\\\hline
$\mathtt{3}$\\\hline
\end{tabular}
\overset{\widetilde{e}_{2}}{\rightarrow}\mathrm{w}%
\begin{tabular}
[c]{|l|}\hline
$\mathtt{1}$\\\hline
$\mathtt{2}$\\\hline
\end{tabular}
.$
\end{example}

In fact $\widetilde{f}_{i_{1}}^{p_{1}},...,\widetilde{f}_{i_{r}}^{p_{r}}$ are
chosen to verify:%
\begin{equation}
f_{i_{1}}^{(p_{1})}\cdot\cdot\cdot f_{i_{r}}^{(p_{r})}v_{\Lambda_{p}%
}=\widetilde{f}_{i_{1}}^{p_{1}}\cdot\cdot\cdot\widetilde{f}_{i_{r}}^{p_{r}%
}v_{\Lambda_{p}}. \label{eq_G_C_for_fund}%
\end{equation}
This implies

\begin{theorem}
\label{Th_marsh}(Marsh) For any admissible column $C$%
\[
G(C)=f_{i_{1}}^{(p_{1})}\cdot\cdot\cdot f_{i_{r}}^{(p_{r})}v_{\Lambda_{p}%
}\text{,}%
\]
where the integers $i_{1},...,i_{r}$ and $p_{1},...,p_{r}$ are determined by
the algorithm above.
\end{theorem}

\subsection{The representation $W(\Lambda_{p})$}

It is well know that the fundamental $sp_{2n}$-module of weight $\Lambda_{p}$
may be regarded as an irreducible component of the $p$-th exterior power of
the vector representation. This exterior power contains a natural basis
indexed by the columns of height $p$.\ Our aim in this subsection is to obtain
a analogous embedding for the $U_{q}(sp_{2n})$-module $V(\Lambda_{p})$.\ We
are going to describe a $U_{q}(sp_{2n})$-module $W(\Lambda_{p})$ containing an
irreducible component isomorphic to $V(\Lambda_{p})$ (that we identify with
$V(\Lambda_{p})$), whose natural basis $\{v_{C}\}$ is indexed by all columns
of $\mathbf{C}(n,p)$. The action of the generators $e_{i},f_{i}$ and
$q^{h_{i}}$ $i=1,..,n$ of $U_{q}(sp_{2n})$ on each vector $v_{C}%
,C\in\mathbf{C}(n,p)$ is easy to describe. So, using Marsh's algorithm, we
will be able to expand the canonical basis of $V(\Lambda_{p})$ on this basis.

Consider the two totally ordered alphabets%
\[
\mathcal{A}_{n}=\{1<\cdot\cdot\cdot<n\}\text{ and }\overline{\mathcal{A}_{n}%
}=\{\overline{n}<\cdot\cdot\cdot<\overline{1}\}\text{.}%
\]
Let $E_{k}^{+}$ be the vector space of dimension $\binom{n}{k}$ with basis
$B_{k}^{+}=\{v_{C}\}$ where $C_{+}$ runs over the set of columns of height $k$
on $\mathcal{A}_{n}$. We define the action of the operators $e_{i},f_{i}$ and
$q^{h_{i}}$ ($i=1,...,n-1)$ on $E_{k}^{+}$ by:%
\begin{align*}
q^{h_{i}}v_{C_{+}}  &  =\left\{
\begin{tabular}
[c]{l}%
$qv_{C_{+}}$ if $i\in C_{+}$ and $i+1\notin C_{+},$\\
$q^{-1}v_{C_{+}}$ if $i+1\in C_{+}$ and $i\notin C_{+},$\\
$v_{C_{+}}$ otherwise.
\end{tabular}
\right. \\
& \\
e_{i}v_{C_{+}}  &  =\left\{
\begin{tabular}
[c]{l}%
$0$ if $i+1\notin C_{+}$ or $i\in C_{+},$\\
$v_{D}$ where $D=C_{+}-\{i+1\}+\{i\}$ otherwise.
\end{tabular}
\right. \\
& \\
f_{i}v_{C_{+}}  &  =\left\{
\begin{tabular}
[c]{l}%
$0$ if $i\notin C_{+}$ or $i+1\in C_{+},$\\
$v_{D_{+}}$ where $D_{+}=C_{+}-\{i\}+\{i+1\}$ otherwise.
\end{tabular}
\right.
\end{align*}
This endows $E_{k}^{+}$ with the structure of a fundamental $U_{q}(sl_{n}%
)$-module which may be regarded as a $q$-analogue of the fundamental $sl_{n}%
$-module $\Lambda^{k}\mathbb{C}$.\ Moreover the basis $\{v_{C}\}$ is then the
canonical basis and the crystal basis of $E_{k}^{+}$ because the fundamental
modules of $U_{q}(sl_{n})$ are minuscule. Similarly $E_{k}^{-}$ the vector
space of dimension $\binom{n}{k}$ with basis $B_{k}^{-}=\{v_{C_{-}}\}$ where
$C_{-}$ runs over the set of columns of height $k$ on $\overline{\mathcal{A}%
}_{n}$ is a $q$-analogue of the fundamental $sl_{n}$-module $\Lambda
^{n-k}\mathbb{C}$ once defined the appropriate action of $e_{i},f_{i}$ and
$q^{h_{i}}$ (obtained by replacing in the right hand sides of the above
formulas $i$ by $\overline{i+1}$ and $i+1$ by $\overline{i}$). Then $B_{k}%
^{-}=\{v_{C_{-}}\}$ is again the canonical basis and the crystal basis of
$E_{k}^{-}$.

Set $\Omega^{+}=\overset{n}{\underset{k=0}{%
{\textstyle\bigoplus}
}}E_{k}^{+}$, $\Omega^{-}=\overset{n}{\underset{k=0}{%
{\textstyle\bigoplus}
}}E_{k}^{-}$ and write $B^{+}=\underset{k=0}{\overset{n}{\cup}}B_{k}^{+}$,
$B^{-}=\underset{k=0}{\overset{n}{\cup}}B_{k}^{-}$.\ Then $B^{+}%
{\textstyle\bigotimes}
B^{-}=\{v_{C^{+}}%
{\textstyle\bigotimes}
v_{C^{-}};$ $v_{C^{+}}\in B^{+},$ $v_{C^{-}}\in B^{-}\}$ is a basis of
$\Omega^{+}%
{\textstyle\bigotimes}
\Omega^{-}$. For $p\in\{1,...,n\},$ let $W(\Lambda_{p})$ be the subspace of
$\Omega^{+}%
{\textstyle\bigotimes}
\Omega^{-}$ generated by the vectors $v_{C_{+}}%
{\textstyle\bigotimes}
v_{C_{-}}$ such that $h(C_{+})+h(C_{-})=p.$ For each of these vectors we set
$v_{C_{+}}%
{\textstyle\bigotimes}
v_{C_{-}}=v_{C}$ where $C$ is the column of height $p$ on $\mathcal{C}_{n}$
(admissible or not) consisting of the letters of $C_{+}$ and $C_{-}$. Then
$\{v_{C}$ $\}$ is a basis of $W(\Lambda_{p}).$ The action of $U_{q}(sl_{n})$
on this basis is deduced from the description of its action on $B_{+}$ and
$B_{-}$ and from (\ref{tensor1}), (\ref{tensor2}) and (\ref{tensor3}). The
formulas below describe this action on the basis vector $v_{C}.$ Set
$E=C\cap\{\overline{i+1},\overline{i},i,i+1\}$ then:%
\begin{equation}
f_{i}v_{C}=\left\{
\begin{tabular}
[c]{l}%
\textrm{(i)}: $v_{D}$ with $D=C-\{i\}+\{i+1\}$ if $E=\{i\},$\\
\\
\textrm{(ii)}: $v_{D}$ with $D=C-\{i\}+\{i+1\}$ if $E=\{\overline
{i+1},\overline{i},i\},$\\
\\
\textrm{(iii)}: $v_{D}$ with $D=C-\{\overline{i+1}\}+\{\overline{i}\}$ if
$E=\{\overline{i+1}\},$\\
\\
\textrm{(iv)}: $v_{D}$ with $D=C-\{\overline{i+1}\}+\{\overline{i}\}$ if
$E=\{\overline{i+1},i,i+1\},$\\
\\
\textrm{(v)}: $q^{-1}v_{D}$ with $D=C-\{\overline{i+1}\}+\{\overline{i}\}$ if
$E=\{\overline{i+1},i+1\},$\\
\\
\textrm{(vi)}: $v_{D}$ with $D=C-\{i\}+\{i+1\}$ if $E=\{\overline{i},i\},$\\
\\
\textrm{(vii)}: $v_{D_{1}}+qv_{D_{2}}$ with$\left\{
\begin{tabular}
[c]{l}%
$D_{1}=C-\{i\}+\{i+1\}$\\
$D_{2}=C-\{\overline{i+1}\}+\{\overline{i}\}$%
\end{tabular}
\right.  $if $E=\{\overline{i+1},i\},$\\
\\
\textrm{(viii)}: $0$ otherwise.
\end{tabular}
\right.  \label{fi}%
\end{equation}%
\begin{equation}
e_{i}v_{C}=\left\{
\begin{tabular}
[c]{l}%
\textrm{(i)}: $v_{D}$ with $D=C-\{i+1\}+\{i\}$ if $E=\{i+1\},$\\
\\
\textrm{(ii)}: $v_{D}$ with $D=C-\{i+1\}+\{i\}$ if $E=\{\overline
{i+1},\overline{i},i+1\},$\\
\\
\textrm{(iii)}: $v_{D}$ with $D=C-\{\overline{i}\}+\{\overline{i+1}\}$ if
$E=\{\overline{i}\},$\\
\\
\textrm{(iv)}: $v_{D}$ with $D=C-\{\overline{i}\}+\{\overline{i+1}\}$ if
$E=\{\overline{i},i,i+1\},$\\
\\
\textrm{(v)}: $q^{-1}v_{D}$ with $D=C-\{i+1\}+\{i\}$ if $E=\{\overline
{i+1},i+1\},$\\
\\
\textrm{(vi)}: $v_{D}$ with $D=C-\{\overline{i}\}+\{\overline{i+1}\}$ if
$E=\{\overline{i},i\},$\\
\\
\textrm{(vii)}: $v_{D_{1}}+qv_{D_{2}}$ with$\left\{
\begin{tabular}
[c]{l}%
$D_{1}=C-\{\overline{i}\}+\{\overline{i+1}\}$\\
$D_{2}=C-\{i+1\}+\{i\}$%
\end{tabular}
\right.  $ if $E=\{\overline{i},i+1\},$\\
\\
\textrm{(viii)}: $0$ otherwise.
\end{tabular}
\right.  \label{ei}%
\end{equation}%

\begin{equation}
q^{h_{i}}v_{C}=q^{<\mathrm{wt}(C),\Lambda_{i}>}\,v_{C}. \label{hi}%
\end{equation}
Now we define an action of the operators $e_{n},f_{n}$ and $q^{h_{n}}$ in
order to endow $W(\Lambda_{p})$ with the structure of a $U_{q}(sp_{2n}%
)$-module. Set%
\begin{gather}
f_{n}v_{C}=\left\{
\begin{tabular}
[c]{l}%
$0$ if $\overline{n}\in C$ or $n\notin C$\\
$v_{D}$ with $D=C-\{n\}+\{\overline{n}\}$ otherwise
\end{tabular}
\right.  ,\label{f0}\\
\nonumber\\
e_{n}v_{C}=\left\{
\begin{tabular}
[c]{l}%
$0$ if $\overline{n}\notin C$ or $n\in C$\\
$v_{D}$ with $D=C-\{\overline{n}\}+\{n\}$ otherwise
\end{tabular}
\right.  ,\label{e0}\\
\nonumber\\
q^{h_{n}}v_{C}=q^{<\mathrm{wt}(C),\Lambda_{n}>}\,v_{C}. \label{h0}%
\end{gather}

\begin{lemma}
The above actions of $e_{n},f_{n}$ and $q^{h_{n}}$ make $W(\Lambda_{p})$ into
a $U_{q}(sp_{2n})$-module.
\end{lemma}

\begin{proof}
It suffices to show that the actions of $f_{n},e_{n}$ and $h_{n}$ are
compatible with (\ref{D1}), (\ref{D2}), (\ref{D3}), (\ref{D4}), (\ref{D5}) and
(\ref{D6}). $W(\Lambda_{p})$ is already a $U_{q}(sl_{n})$-module so we can
restrict ourselves to the cases where $i=n$ or $j=n$ in these relations.
Suppose first $(i,j)=(n,n)$. Then denote by $U_{q}^{0}(sl_{2})$ the subalgebra
of $U_{q}(sp_{2n})$ generated by $f_{n},e_{n}$ and $q^{h_{n}}$. Formulas
(\ref{f0}),(\ref{e0}) and (\ref{h0}) imply that $W(\Lambda_{p})$ is a
$U_{q}^{0}(sl_{2})$-module and the only non trivial module occurring in its
decomposition into irreducible $U_{q}^{0}(sl_{2})$-modules is the vector
representation of $U_{q}^{0}(sl_{2})$. Hence we can suppose $(i,j)\neq(n,n)$.
Let $C\in\mathbf{C}(n,p)$. For $i\neq n$, we have to establish that:%
\begin{gather*}
\text{\textrm{(i)}}:\text{ }q^{h_{n}}e_{i}q^{-h_{n}}v_{C}=\left\{
\begin{tabular}
[c]{l}%
$q^{-1}e_{i}v_{C}$ if $i=n-1$\\
$e_{i}v_{C}$ otherwise
\end{tabular}
\right.  \text{ and }q^{h_{n}}f_{i}q^{-h_{n}}v_{C}=\left\{
\begin{tabular}
[c]{l}%
$qe_{i}v_{C}$ if $i=n-1$\\
$f_{n}v_{C}$ otherwise
\end{tabular}
\right.  ,\\
\text{\textrm{(ii)}}:\text{ }q^{h_{i}}e_{n}q^{-h_{i}}v_{C}=\left\{
\begin{tabular}
[c]{l}%
$q^{-2}e_{n}v_{C}$ if $i=n-1$\\
$e_{n}v_{C}$ otherwise
\end{tabular}
\right.  \text{ and }q^{h_{i}}f_{n}q^{-h_{i}}v_{C}=\left\{
\begin{tabular}
[c]{l}%
$q^{2}e_{n}v_{C}$ if $i=n-1$\\
$f_{n}v_{C}$ otherwise
\end{tabular}
\right.  ,\\
\text{\textrm{(iii)}}:\text{ }[e_{n},f_{i}]v_{C}=[f_{n},e_{i}]v_{C}=0\\
\mathrm{(iv)}:\text{ }[e_{n},e_{i}]v_{C}=[f_{n},f_{i}]v_{C}=0\text{ for }i\neq
n-1,\\
\text{\textrm{(v)}}:\text{{}}\left\{
\begin{tabular}
[c]{l}%
$(e_{n-1}e_{n}^{2}-(q^{2}+q^{-2})e_{n}e_{n-1}e_{n}+e_{n}^{2}e_{n-1})v_{C}=0$\\
$(f_{n-1}f_{n}^{2}-(q^{2}+q^{-2})f_{n}f_{n-1}f_{n}+f_{n}^{2}e_{n-1})v_{C}=0$%
\end{tabular}
\right.  ,\\
\text{\textrm{(vi)}}:\text{ }\left\{  \text{%
\begin{tabular}
[c]{l}%
$\lbrack e_{n}e_{n-1}^{3}-(q^{2}+1+q^{-2})e_{n-1}e_{n}e_{n-1}^{2}%
+(q^{2}+1+q^{-2})e_{n-1}^{2}e_{n}e_{n-1}-e_{n-1}^{3}e_{n}]v_{C}=0$\\
$\lbrack f_{n}f_{n-1}^{3}-(q^{2}+1+q^{-2})f_{n-1}f_{n}f_{n-1}^{2}%
+(q^{2}+1+q^{-2})f_{n-1}^{2}f_{n}f_{n-1}-f_{n-1}^{3}f_{n}]v_{C}=0$%
\end{tabular}
}\right.  .
\end{gather*}
When $i\neq n-1$, the actions the operators $e_{n},f_{n}$ or $q^{h_{n}}$ on
$W(\Lambda_{p})$ commute with those of the operators $e_{i},f_{i}$ or
$q^{h_{i}}.$ Moreover we have $e_{i}^{3}v_{C}=f_{i}^{3}v_{C}=0$ and $e_{n}%
^{2}v_{C}=f_{n}^{2}v_{C}=0$. This implies that the above relations are
immediate when $i\neq n-1$. If $i=n-1,$ we obtain \textrm{(i)} and
\textrm{(ii)} by a case by case computation from (\ref{ei}), (\ref{fi}) and
(\ref{h0}). Relations \textrm{(iii)} follows from the equalities $e_{n-1}%
f_{n}v_{C}=f_{n}e_{n-1}v_{C}=e_{n}f_{n-1}v_{C}=f_{n-1}e_{n}v_{C}=0$. To
establish relations \textrm{(v)} and \textrm{(vi)} it suffices to prove that%
\begin{gather*}
e_{n}e_{n-1}e_{n}v_{C}=f_{n}f_{n-1}f_{n}v_{C}=0\\
f_{n-1}f_{n}f_{n-1}^{2}v_{C}=f_{n-1}^{2}f_{n}f_{n-1}v_{C}=0;\text{
\ \ \ \ \ \ \ \ \ }e_{n-1}e_{n}e_{n-1}^{2}v_{C}=e_{n-1}^{2}e_{n}e_{n-1}%
v_{C}=0.
\end{gather*}
If $f_{n}v_{C}=v_{D}\neq0$, then the column $D$ contains $\overline{n}$ but
not $n$ so $D$ is of type \textrm{(ii)}, \textrm{(iii)} or \textrm{(vii)} in
\ref{fi}$.$ Then we obtain after a simple computation that $f_{n}f_{n-1}%
f_{n}v_{C}=0.$ We prove similarly the equality: $e_{n}e_{n-1}e_{n}v_{C}=0$
which implies \textrm{(v)}. From \textrm{(vii)} of formula (\ref{fi}) it
follows that $f_{n}f_{n-1}^{2}v_{C}=0$. Moreover $f_{n-1}^{2}f_{n}f_{n-1}%
v_{C}=0$ because $f_{n}f_{n-1}v_{C}=0$ or may be written $v_{D}$ with $D$ a
column which is not of the type \textrm{(vii)} of (\ref{fi}). In the same
manner we can see that $e_{n}e_{n-1}^{2}v_{C}=0=e_{n-1}^{2}e_{n}e_{n-1}v_{C}$
which implies \textrm{(vi)}.
\end{proof}

$W(\Lambda_{p})$ is not an irreducible $U_{q}(sp_{2n})$-module. Since it is
finite-dimensional, it decomposes into a direct sum of irreducible components.
It is easy to verify that the vector $v_{C_{p}^{0}}$ where $C_{p}^{0}$ is the
column defined in \ref{subsec_crystal} is a highest weight vector of weight
$\Lambda_{p}$. Hence setting $v_{\Lambda_{p}}=v_{C_{p}^{0}},$ we can identify
$U_{q}(sp_{2n})v_{C_{p}^{0}}$ and $V(\Lambda_{p})$. Then it is possible to
expand explicitly the vectors of the canonical basis of $V(\Lambda_{p})$ on
the basis $\{v_{C};C\in\mathbf{C}(n,p)\}$.

\begin{lemma}
\ \ \ \ \ \ \ \ \label{Lem_base_cryst_Vp}

\begin{enumerate}
\item  Let $\mathcal{L}_{p}$ be the $A$-submodule of $W(\Lambda_{p})$
generated by the vectors $v_{C},$ $C\in\mathbf{C}(n,p)$. Write $\mathrm{w(}C)$
for the image of the vector $v_{C},$ $C\in\mathbf{C}(n,p)$ by the projection
$\pi_{p}:\mathcal{L}_{p}\rightarrow\mathcal{L}_{p}/q\mathcal{L}_{p}$ and
$\mathcal{B}_{p}=\{\mathrm{w(}C)\left|  {}\right.  C\in\mathbf{C}(n,p)\}$.
Then $(\mathcal{L}_{p},\mathcal{B}_{p})$ is a crystal basis of $W(\Lambda
_{p})$.

\item  Set $L(\Lambda_{p})=\mathcal{L}_{p}\cap V(\Lambda_{p})$ and denote by
$B(\Lambda_{p})$ the set of images of the vectors $\{v_{C};$ $C\in
\mathbf{Ca}(n,p)\}$ by $\pi_{p}$. Then $(L(\Lambda_{p}),B(\Lambda_{p}))$ is
the crystal basis of $V(\Lambda_{p})$.
\end{enumerate}
\end{lemma}

\begin{proof}
$1:$ We have to prove that $\mathcal{L}_{p}$ and $\mathcal{B}_{p}$ verify the
assertions $\mathrm{(i)},$\textrm{ }$\mathrm{(ii)},$\textrm{ }$\mathrm{(iii)}$
and $\mathrm{(iv)}$ of (\ref{def_C_B}).\ When $i=n$ this follows immediately
from formulas (\ref{f0}) and (\ref{e0}). Indeed the actions of $\widetilde
{f}_{n}$ and $f_{n}$ (resp. $\widetilde{e}_{n}$ and $e_{n}$) coincide on each
vector $v_{C}$, $C\in\mathbf{C}(n,p).\;$On the other hand, by construction
$W(\Lambda_{p})$ is a direct sum of tensor products of $U_{q}(sl_{n})$-modules
and it follows from the compatibility of crystals with tensor products that
$(\mathcal{L}_{p},\mathcal{B}_{p})$ is a crystal basis of $W(\Lambda_{p})$
considered as a $U_{q}(sl_{n})$-module.\ Hence the assertions $\mathrm{(i)}%
,$\textrm{ }$\mathrm{(ii)},$\textrm{ }$\mathrm{(iii)}$ and $\mathrm{(iv)}$ are
also verified for $i=1,...,n-1$. Note that the actions of $\widetilde{e}_{i}$
and $\widetilde{f}_{i},$ $i=1,...,n$ on $\mathcal{B}_{p}$ coincide with those
given by Rule \ref{+-}.

$2:$ It is clear that $L(\Lambda_{p})$ verifies the conditions of
(\ref{def_C_B}).$\;$By theorem 4.2 of \cite{Ka2}, we know that
\[
B(\Lambda_{p})=\{\widetilde{f}_{i_{1}}^{a_{1}}\cdot\cdot\cdot\widetilde
{f}_{i_{r}}^{a_{r}}\mathrm{w}(C_{p}^{0});\text{ }i_{1},...,i_{r}%
=1,...,n;\text{ }a_{1},...,a_{r}>0\}-\{0\}
\]
for $\mathrm{w}(C_{p}^{0})$ is the image of $v_{\Lambda_{p}}$ by $\pi_{p}$. We
have just seen that the actions of $\widetilde{e}_{i}$ and $\widetilde{f}%
_{i},$ $i=1,...,n$ on $\mathcal{B}_{p}$ coincide with those given by Rule
\ref{+-}.$\;$Hence, by Theorem \ref{TH_KN}, $B(\Lambda_{p})=\{\mathrm{w}(C);$
$C\in\mathbf{Ca}(n,p)\}$.
\end{proof}

\subsection{Expression of $\{G(C)\}$ on $\{v_{C}\mathbf{\}}$}

Let $C$ be an admissible column of height $p$ and set $C=\widetilde{f}_{i_{1}%
}^{p_{1}}\cdot\cdot\cdot\widetilde{f}_{i_{r}}^{p_{r}}C_{p}^{0}$ with the
notation of Section \ref{subsec_marsh}. Then by Theorem \ref{Th_marsh}%
\[
G(C)=f_{i_{1}}^{(p_{1})}\cdot\cdot\cdot f_{i_{r}}^{(p_{r})}v_{\Lambda_{p}%
}\text{.}%
\]
We are going to give a combinatorial description of $G(C)$. Define the
$r$-tuple $L_{C}=(u_{1},...,u_{r})\subset\mathcal{A}_{n}$ from the $r$-tuple
$K_{C}=(x_{1}<\cdot\cdot\cdot<x_{r})$ containing the unbarred letters $x$ such
that $(x,\overline{x})\subset C$ by:%
\begin{gather*}
u_{1}<x_{1},\{u_{1},\overline{u}_{1}\}\cap C=\text{\O\ }\ \text{and }%
u_{1}\text{ is maximal,}\\
\text{for }i=2,...,r,\text{ }u_{i}<x_{i},\text{ }\{u_{i},\overline{u}%
_{i}\}\cap(C\cup\{u_{1},...,u_{i-1}\})=\text{\O\ \ and }u_{i}\text{ is
maximal.}%
\end{gather*}
Then the letters of $L_{C\,}$ are those of the set $J_{C}$ defined in
\ref{subsec_crystal}. Indeed we have seen that $J_{C}$ may be regarded as the
set of $r$ maximal unbarred letters $\{t_{1},...,t_{r}\}$ such that
$t_{i}<z_{i}$ and $\{t_{i},\overline{t}_{i}\}\cap C=\emptyset$.\ Notice that,
in general, the letters of $L_{C}$ are not ordered in decreasing order as
those of $J_{C}$.

\begin{example}
$C=%
\begin{tabular}
[c]{|c|}\hline
$\mathtt{3}$\\\hline
$\mathtt{5}$\\\hline
$\mathtt{6}$\\\hline
$\mathtt{\bar{6}}$\\\hline
$\mathtt{\bar{5}}$\\\hline
$\mathtt{\bar{3}}$\\\hline
\end{tabular}
$ is admissible with $K_{C}=(3<5<6)$ and $L_{C}=(2,4,1).$
\end{example}

For any subset $X=\{x_{i_{1}},...,x_{i_{s}}\}\subset K_{C},$ let $C_{X}$ be
the column of $\mathbf{C}(n,p)$ obtained by changing in $C$ each pair of
letters $(x_{i_{j}},$ $\overline{x_{i_{j}}})$ into the corresponding pair of
letters $(u_{i_{j}},\overline{u_{i_{j}}})$. Then, with the above notations we obtain:

\begin{theorem}
\label{Th_G(C)}For any admissible column $C$ of height $p$%
\[
G(C)=\underset{X\subset K_{C}}{\sum}q^{\mathrm{card}(X)}v_{C_{X}}.
\]
\end{theorem}

\begin{example}
With the column $C$ of the previous example we have:\vspace{0.5cm}%

\[
G(C)=%
\begin{tabular}
[c]{|c|}\hline
$\mathtt{3}$\\\hline
$\mathtt{5}$\\\hline
$\mathtt{6}$\\\hline
$\mathtt{\bar{6}}$\\\hline
$\mathtt{\bar{5}}$\\\hline
$\mathtt{\bar{3}}$\\\hline
\end{tabular}
+q\left(  \text{ }%
\begin{tabular}
[c]{|c|}\hline
$\mathtt{2}$\\\hline
$\mathtt{5}$\\\hline
$\mathtt{6}$\\\hline
$\mathtt{\bar{6}}$\\\hline
$\mathtt{\bar{5}}$\\\hline
$\mathtt{\bar{2}}$\\\hline
\end{tabular}
+%
\begin{tabular}
[c]{|c|}\hline
$\mathtt{3}$\\\hline
$\mathtt{4}$\\\hline
$\mathtt{6}$\\\hline
$\mathtt{\bar{6}}$\\\hline
$\mathtt{\bar{4}}$\\\hline
$\mathtt{\bar{3}}$\\\hline
\end{tabular}
+%
\begin{tabular}
[c]{|c|}\hline
$\mathtt{1}$\\\hline
$\mathtt{3}$\\\hline
$\mathtt{5}$\\\hline
$\mathtt{\bar{5}}$\\\hline
$\mathtt{\bar{3}}$\\\hline
$\mathtt{\bar{1}}$\\\hline
\end{tabular}
\text{ }\right)  +q^{2}\left(  \text{ }%
\begin{tabular}
[c]{|c|}\hline
$\mathtt{2}$\\\hline
$\mathtt{4}$\\\hline
$\mathtt{6}$\\\hline
$\mathtt{\bar{6}}$\\\hline
$\mathtt{\bar{4}}$\\\hline
$\mathtt{\bar{2}}$\\\hline
\end{tabular}
+%
\begin{tabular}
[c]{|c|}\hline
$\mathtt{1}$\\\hline
$\mathtt{2}$\\\hline
$\mathtt{5}$\\\hline
$\mathtt{\bar{5}}$\\\hline
$\mathtt{\bar{2}}$\\\hline
$\mathtt{\bar{1}}$\\\hline
\end{tabular}
+%
\begin{tabular}
[c]{|c|}\hline
$\mathtt{1}$\\\hline
$\mathtt{4}$\\\hline
$\mathtt{3}$\\\hline
$\mathtt{\bar{3}}$\\\hline
$\mathtt{\bar{4}}$\\\hline
$\mathtt{\bar{1}}$\\\hline
\end{tabular}
\text{ }\right)  +q^{3}\,%
\begin{tabular}
[c]{|c|}\hline
$\mathtt{1}$\\\hline
$\mathtt{2}$\\\hline
$\mathtt{4}$\\\hline
$\mathtt{\bar{4}}$\\\hline
$\mathtt{\bar{2}}$\\\hline
$\mathtt{\bar{1}}$\\\hline
\end{tabular}
\]
where we have written for short $C$ in place of $v_{C}$.
\end{example}

\bigskip

\begin{proof}
Let $C$ be an admissible column and set $\mathrm{w(}C)=\widetilde{f}_{i_{1}%
}^{p_{1}}\cdot\cdot\cdot\widetilde{f}_{i_{r}}^{p_{r}}\mathrm{w(}C_{p}^{0})$
with the notation of Section \ref{subsec_marsh}. First notice that by
(\ref{fi}) and (\ref{f0}) we will have $d(D)=d(C)$ for any column $D$ such
that $v_{D}$ occurs with a non-zero coefficient in the decomposition of $G(C)$
on the basis $\{v_{C}\}$.\ We proceed by induction on $r$. The theorem is
clear for $r=0,$ i.e. $C=C_{p}^{0}$.\ Suppose the result is proved for the
admissible column $\mathrm{w(}C^{\prime})=\widetilde{f}_{i_{2}}^{p_{2}}%
\cdot\cdot\cdot\widetilde{f}_{i_{r}}^{p_{r}}\mathrm{w(}C_{p}^{0}),$ that is:%
\[
G(C^{\prime})=\underset{X^{\prime}\subset K_{C^{\prime}}}{\sum}%
q^{\mathrm{card}(X^{\prime})}v_{C_{X^{\prime}}^{\prime}}.
\]
When $p_{1}=2,$ the letters $i_{1}$ and $\overline{i_{1}+1}$ occur in
$C^{\prime}$ but not the letters $\overline{i}_{1}$ or $i_{1}+1$. We have the
same property for all the columns $C_{X^{\prime}}^{\prime}$. Moreover
$C=C^{\prime}-\{i,\overline{i+1}\}+\{i+1,\overline{i}\},$ $K_{C}=K_{C^{\prime
}}$ and $L_{C}=L_{C^{\prime}}$.\ Then $G(C)=f_{i_{1}}^{(2)}G(C^{\prime})$ is
obtained by replacing in the columns $C_{X^{\prime}}^{\prime}$ each pair
$(i,\overline{i+1})$ by $(i+1,\overline{i})$. So the theorem is proved.

When $p_{1}=1$ and $i_{1}=n$, $n$ $\in C_{X^{\prime}}^{\prime}$ but
$\overline{n}\notin C_{X^{\prime}}^{\prime}$. Then $K_{C}=K_{C^{\prime}}$,
$L_{C}=L_{C^{\prime}}$ and $G(C)$ is obtained by replacing all the letters $n$
by letters $\overline{n}$ into the columns $C_{X^{\prime}}^{\prime}$. So the
theorem is proved.

When $p_{1}=1$ and $i_{1}\neq n$, $C$ can not contain the pair $(i+1,\overline
{i})$ without containing a letter of $\{i_{1},\overline{i_{1}+1}\}$ (otherwise
$p_{1}=2$). Hence $C^{\prime}$ may not be of type \textrm{(v)} in (\ref{fi}).
Moreover $C^{\prime}$ is not of type \textrm{(vi)} or \textrm{(viii)
}otherwise $\widetilde{f}_{i_{1}}\mathrm{w(}C^{\prime})=0$. Suppose
$C^{\prime}$ of type \textrm{(i) }or \textrm{(ii)}. Then $C=C^{\prime}%
-\{i_{1}\}+\{i_{1}+1\},$ $K_{C}=K_{C^{\prime}}$ and $L_{C}=L_{C^{\prime}}%
$.\ Hence
\[
G(C)=f_{i_{1}}(G(C^{\prime}))=\underset{X^{\prime}\subset K_{C}}{\sum
}q^{\mathrm{card}(X)}v_{C_{X^{\prime}}^{^{\prime}}-\{i_{1}\}+\{i_{1}%
+1\}}=\underset{X\subset K_{C}}{\sum}q^{\mathrm{card}(X)}v_{C_{X}}.
\]
When $C^{\prime}$ of type \textrm{(iii) }or \textrm{(iv) }we have\textrm{
}$C=C^{\prime}-\{\overline{i_{1}+1}\}+\{\overline{i_{1}}\},$ $K_{C}%
=K_{C^{\prime}}$ and $L_{C}=L_{C^{\prime}}$.\ Hence%
\[
G(C)=f_{i_{1}}(G(C^{\prime}))=\underset{X^{\prime}\subset K_{C}}{\sum
}q^{\mathrm{card}(X)}v_{C_{X^{\prime}}^{\prime}-\{\overline{i_{1}%
+1}\}+\{\overline{i_{1}}\}}=\underset{X\subset K_{C}}{\sum}q^{\mathrm{card}%
(X)}v_{C_{X}}.
\]
If $C^{\prime}$ is of type \textrm{(vii)} in (\ref{fi}), the letters $i_{1}$
and $\overline{i_{1}+1}$ occur in all the columns $C_{X^{\prime}}^{\prime}$
but not the letters $\overline{i}_{1}$ or $i_{1}+1$. Then $C=C^{\prime
}-\{i_{1}\}+\{i_{1}+1\},$ $K_{C}=K_{C^{\prime}}+\{i_{1}+1\}$. Note that
$i_{1}+1$ is the lowest letter of $K_{C}$ because $i_{1}+1$ is the lowest
movable letter of $C$. Hence $L_{C}=L_{C^{\prime}}+\{i_{1}\}$ and the letter
of $L_{C}$ corresponding to $i_{1}+1\in K_{C}$ is $i_{1}$.\ We obtain
$f_{i_{1}}v_{C_{X^{\prime}}^{\prime}}=v_{C_{X^{\prime}}}+qv_{C_{X^{\prime
}+\{i_{1}+1\}}}$ which implies:%
\[
G(C)=\underset{X^{\prime}\subset K_{C^{\prime}}}{\sum}q^{\mathrm{card}%
(X^{\prime})}(v_{C_{X^{\prime}}}+qv_{C_{X^{\prime}+\{i_{1}+1\}}}%
)=\underset{X\subset K_{C}}{\sum}q^{\mathrm{card}(X)}v_{C_{X}}%
\]
because the parts of $K_{C}$ are exactly the elements of the set $\{X^{\prime
},X^{\prime}+\{i_{1}+1\};$ $X^{\prime}\subset K_{C^{\prime}}\}$.
\end{proof}

\section{Computation of the canonical basis of $V(\lambda)\label{sec_G_lambda}%
$}

\subsection{The representation $W(\lambda)$}

Let $\lambda=\underset{p=1}{\overset{n}{\sum}}\lambda_{p}\Lambda_{p}$ be a
dominant weight and write%
\[
W(\lambda)=W(\Lambda_{1})^{\bigotimes\lambda_{1}}%
{\textstyle\bigotimes}
\cdot\cdot\cdot%
{\textstyle\bigotimes}
W(\Lambda_{n})^{\bigotimes\lambda_{n}}.
\]
The natural basis of $W(\lambda)$ consists of the tensor products $v_{C_{r}}%
{\textstyle\bigotimes}
\cdot\cdot\cdot%
{\textstyle\bigotimes}
v_{C_{1}}$ of basis vectors $v_{C}$ of the previous section.\ The
juxtaposition of the columns $C_{1},...,C_{r}$ is called a tabloid of shape
$\lambda$.\ We can regard it as a filling $\tau$ of the Young diagram of shape
$\lambda,$ the $i$-th column of which is equal to $C_{i}$. We shall write
$v_{\tau}=v_{C_{r}}%
{\textstyle\bigotimes}
\cdot\cdot\cdot%
{\textstyle\bigotimes}
v_{C_{1}}.$ Note that the columns of a tabloid are not necessarily admissible
and there is no condition on the rows. Write $\mathbf{T}(n,\lambda)$ for the
set of tabloids of shape $\lambda$. The reading of the tabloid $\tau
=C_{1}\cdot\cdot\cdot C_{r}\in\mathbf{T}(n,\lambda)$ is the word
\textrm{w(}$\tau)=\mathrm{w}(C_{r})\cdot\cdot\cdot\mathrm{w}(C_{1}%
)\in\mathcal{C}_{n}^{\ast}$. The weight of $\tau$ is the weight $v_{\tau}$.

Let $\mathcal{L}_{\lambda}$ be the $A$-submodule of $W(\lambda)$ generated by
the vectors $v_{\tau},$ $\tau\in\mathbf{T}(n,\lambda)$. We identify the image
of the vector $v_{\tau}$ by the projection $\pi_{\lambda}:\mathcal{L}%
_{\lambda}\rightarrow\mathcal{L}_{\lambda}/q\mathcal{L}_{\lambda}$ with the
word \textrm{w(}$\tau)$. The pair $(\mathcal{L}_{\lambda},\mathcal{B}%
_{\lambda}=\{\mathrm{w(}\tau)\left|  {}\right.  \tau\in\mathbf{T}%
(n,\lambda)\})$ is then a crystal basis of $W(\lambda)$. Indeed by Lemma
\ref{Lem_base_cryst_Vp}, it is the tensor product of the crystal bases of the
representations $W(\Lambda_{p})$ occurring in $W(\lambda)$. Denote by
$v_{\lambda}$ the tensor product in $W(\lambda)$ of the highest weight vectors
of each $W(\Lambda_{p})$. We identify $V(\lambda)$ with the submodule of
$W(\lambda)$ of highest weight vector $v_{\lambda}$. Then, with the above
notations, $v_{\lambda}=v_{T_{\lambda}}$ where $T_{\lambda}$ is the symplectic
tableau of shape $\lambda$ whose $k$-th row is filled by letters $k$. By
Theorem 4.2 of \cite{Ka2}, we know that
\[
B(\lambda)=\{\widetilde{f}_{i_{1}}^{a_{1}}\cdot\cdot\cdot\widetilde{f}_{i_{r}%
}^{a_{r}}\mathrm{w}(T_{\lambda});\text{ }i_{1},...,i_{r}=1,...,n;\text{ }%
a_{1},...,a_{r}>0\}-\{0\}.
\]
The actions of $\widetilde{e}_{i}$ and $\widetilde{f}_{i},$ $i=1,...,n$ on
$\mathcal{B}_{\lambda}$ coincide with those given by Rule \ref{+-} because it
is true on each $\mathcal{B}_{p}$, $p=1,...,n$.$\;$Hence, by Theorem
\ref{TH_KN}, $B(\lambda)=\{\mathrm{w}(T);$ $T\in\mathbf{ST}(n,\lambda)\}$. For
each vector $G(T)$ of the canonical basis of $V(\lambda)$ we will have:%
\[
G(T)\equiv v_{T}\operatorname{mod}q\mathcal{L}_{\lambda}.
\]
The aim of this section is to describe an algorithm computing the
decomposition of the canonical basis $\{G(T);T\in\mathbf{ST}(n,\lambda)\}$
onto the basis $\{v_{\tau};\tau\in\mathbf{T}(n,\lambda)\}$ of $W(\lambda)$.

Let $w_{1}=x_{1}\cdot\cdot\cdot x_{l}$ and $w_{2}=y_{1}\cdot\cdot\cdot y_{l}$
be two distinct words on $\mathcal{C}_{n}$ with the same length and $k$ the
lowest integer such that $x_{k}\neq y_{k}.$ Write $w_{1}\trianglelefteq w_{2}$
if $x_{k}\leq y_{k}$ in $\mathcal{C}_{n}$ and $w_{1}\vartriangleright w_{2}$
else that is, $\trianglelefteq$ is the lexicographic order. Then we endow the
set $\mathbf{T}(n,\lambda)$ with a total ordering by setting:%
\[
\tau_{1}\trianglelefteq\tau_{2}\Longleftrightarrow\mathrm{w(}\tau
_{1})\trianglelefteq\mathrm{w(}\tau_{2}).
\]
Notice that for any tabloid $\tau\in\mathbf{T}(n,\lambda)$, we have
$T_{\lambda}\trianglelefteq\tau.$

We are going to compute the canonical basis $\{G(T);T\in\mathbf{ST}%
(n,\lambda)\}$ in two steps. First we obtain an intermediate basis
$\{A(T);T\in\mathbf{ST}(n,\lambda)\}$ which is fixed by the involution
$\overline{\text{%
\begin{tabular}
[c]{l}%
\ \
\end{tabular}
}}$ (condition (\ref{cond_invo})). Next we correct it in order to have
condition (\ref{cond_cong}). This second step is easy because we can prove
that the transition matrix from $\{A(T);T\in\mathbf{ST}(n,\lambda)\}$ to
$\{G(T);T\in\mathbf{ST}(n,\lambda)\}$ is unitriangular once the symplectic
tableaux and the tabloids are ordered by $\trianglelefteq$.

We start with a general Lemma analogous to Lemma 4.1 of \cite{L-T}. Let $\nu$
be a tabloid and $i,m$ two integers such that $f_{i}^{m}v_{\nu}\neq0$. Suppose
that the vector $v_{\tau}$ appears in the decomposition of $f_{i}^{m}v_{\nu}$
on the basis $\{v_{\tau};$ $\tau\in\mathbf{T}(n,\lambda)\}$ with a non zero
coefficient $\kappa_{\tau}$. Then the tabloid $\tau$ is obtained from $\nu$ by
changing $m$ occurrences of letters of $\{\overline{i+1},i\}$ into the
corresponding letters of $\{\overline{i},i+1\}$. Set \textrm{w(}$\nu
)=x_{1}\cdot\cdot\cdot x_{l}$ and denote by $(x_{i_{1}},...,x_{i_{m}}),$
$i_{1}<\cdot\cdot\cdot<i_{m}$ the $m$-tuple of letters of \textrm{w(}$\nu)$
modified to obtain \textrm{w(}$\tau)$. By formula (\ref{tensor3}), the vector
$v_{\tau}$ will appear $m!$ times in $f_{i}^{m}v_{\nu}$ according to the $m!$
permutation of $(i_{1},...,i_{m}).$ For each permutation $(i_{\sigma
(1)},...,i_{\sigma(m)})$ of $(i_{1},...,i_{m})$ write $q_{i}^{N(\sigma)}$ for
the power of $q_{i}$ appearing when we compute \textrm{w(}$\tau)$ by changing
the letters $x_{i_{\sigma(1)}},...,x_{i_{\sigma(m-1),}}x_{i_{\sigma(m)}}$ of
\textrm{w(}$\nu\mathrm{)}$ in this order. Then we obtain by induction on the
length $l(\sigma)$ of the permutation $\sigma$:%
\[
N(\sigma)=N(\mathrm{id})+2l(\sigma).
\]
So the coefficient of $v_{\tau}$ in $f_{i}^{s}v_{\nu}$ is equal to
\begin{equation}
\kappa_{\tau}=\underset{\sigma}{\sum}q_{i}^{N(\mathrm{id})+2l(\sigma)}%
=q_{i}^{N(\mathrm{id})+\frac{m(m+1)}{2}}[m]_{i}. \label{coeff_N_sigma}%
\end{equation}
Hence the coordinates of $f_{i}^{(m)}v_{\nu}=f_{i}^{m}v_{\nu}/[m]_{i}$ belongs
to $\mathbb{N}[q,q^{-1}]$. The following lemma follows by induction on $s$:

\begin{lemma}
\label{Lem_F()}Let $v\in V(\lambda)$ be a vector of the type
\begin{equation}
v=f_{i_{1}}^{(r_{1})}\cdot\cdot\cdot f_{i_{s}}^{(r_{s})}v_{\lambda}
\label{Type_F()}%
\end{equation}
where $(i_{1},...,i_{s})$ and $(r_{1},...,r_{s})$ are two sequences of
integers.\ Then the coordinates of $v$ on the basis $\{v_{\tau};\tau
\in\mathbf{T}(n,\lambda)\}$ belong to $\mathbb{N}[q,q^{-1}]$.
\end{lemma}

\subsection{The basis $\{A(T)\}$}

The basis $\{A(T)\}$ will be a monomial basis, that is, a basis of the form
\begin{equation}
A(T)=f_{i_{1}}^{(r_{1})}\cdot\cdot\cdot f_{i_{s}}^{(r_{s})}v_{\lambda}.
\label{A(T)_monom}%
\end{equation}
By Lemma \ref{Lem_F()}, the coordinates of $A(T)$ on the basis $\{v_{\tau}\}$
of $W(\lambda)$ belong to $\mathbb{N[}q,q^{-1}]$. Let $T=C_{1}\cdot\cdot\cdot
C_{l}\neq T_{\lambda}\in\mathbf{ST}(n,\lambda).$ To find the two sequences of
integers $(i_{1},...,i_{s})$ and $(r_{1},...,r_{s})$ associated to $T$, we
proceed as follows. Let $C_{k}$ be the rightmost column of $T$ such that
$C_{k}$ is not of highest weight, $x$ the lowest movable letter of $C_{k}$ and
$i\in\{1,...,n\}$ such that $\widetilde{e}_{i}(x)\neq0$. Denote by $y$ the
rightmost letter of \textrm{w(}$T)$ such that $y\in\{\overline{i},i+1\}$ and
the factor $x\cdot\cdot\cdot y$ of \textrm{w(}$T)$ contains no letter of
$\{\overline{i+1},i\}$. Set $r_{1}$ for the number of letters of
$\{\overline{i},i+1\}$ in $x\cdot\cdot\cdot y.$ Then $i_{1}$ is defined to be
$i$ and $T_{1}$ is defined to be the tabloid obtained by changing in
$x\cdot\cdot\cdot y$ each letter $\overline{i}$ into $\overline{i+1},$ and
each letter $i+1$ into $i$ if $i\neq n$, and each letter $\overline{n}$ into
$n$ if $i=n$.

\begin{lemma}
\label{T1_Tab}$T_{1}\in\mathbf{ST}(n,\lambda).$
\end{lemma}

Then we do the same with $T_{1}$ getting a new symplectic tableau $T_{2}$ and
a new integer $i_{2}$. And so on until the tableau $T_{s}$ obtained is equal
to $T_{\lambda}$. Notice that we can not write \textrm{w(}$T_{1}%
)=\widetilde{e}_{i_{1}}^{r_{1}}\mathrm{w(}T)$ in general, that is, our
algorithm does not provide a path in the crystal graph $B(\lambda)$ joining
the vertex $\mathrm{w(}T)$ to the vertex of highest weight $\mathrm{w(}%
T_{\lambda})$.

\bigskip

\begin{proof}
(of Lemma \ref{T1_Tab}).

Set $T_{1}=D_{1}\cdot\cdot\cdot D_{l}.$ If $T$ does not contain a letter of
$\{\overline{i+1},i\}$ we may write \textrm{w}$(T_{1})=\widetilde{e}%
_{i}^{r_{1}}\mathrm{w}(T)$ by Rule \ref{+-}.\ So the lemma is true in this
case. Otherwise, let $C_{m}$ be the rightmost column of $T$ containing a
letter of $\{\overline{i+1},i\}.$

The letter $x\in\{\overline{i},i+1\}$ is movable in $C_{k}$ so there is no
letter $\overline{i+1}$ or $i$ to the left of $x$ in \textrm{w}$(C_{k}%
)$.\ Indeed if $x=\overline{i}$, then $i\notin C_{k}$ by Remark \ref{rem_util}%
. There is no letter $\overline{i+1}$ or $i$ to the left of $x$ in
\textrm{w}$(T)$. Indeed the columns to the right of $C_{k}$ are of highest
weight so can not contain the barred letter $\overline{i+1}$. Moreover if $i$
appears in \textrm{w}$(T)$ to the left of $x$ in a column word \textrm{w}%
$(C^{\prime})$, the letters $1,2,...,i\in C_{k}$ because $T$ is a symplectic
tableau. So the letter $i+1$ can not be movable in $C_{k}$ and if
$\overline{i}\in C_{k}$ the column $C_{k}$ is not admissible. Set $T^{\prime
}=C_{m+1}\cdot\cdot\cdot C_{l}$ and $T_{1}^{\prime}=D_{m+1}\cdot\cdot\cdot
D_{l}$.\ Rule \ref{+-} implies that \textrm{w}$(T_{1}^{\prime})=\widetilde
{e}_{i}^{p_{1}}\mathrm{w}(T^{\prime})$ with $p_{1}\leq r_{1}$. Hence
$T_{1}^{\prime}$ is a symplectic tableau.

Suppose first that $C_{m}=D_{m}$ (hence $y\notin C_{m}$ and $p_{1}=r_{1}$) .
Then the columns $C_{1},\cdot\cdot\cdot,C_{m}$ are not modified when $T_{1}$
is computed from $T$. So it suffices to prove that $rC_{m}\leq lD_{m+1}$. If
$\widetilde{e}_{i}\mathrm{w}(C_{m})=0,$ $\widetilde{e}_{i}^{r_{1}}%
(\mathrm{w}(C_{m}T^{\prime}))=\widetilde{e}_{i}^{r_{1}}\mathrm{w}%
(T_{1}^{\prime})\mathrm{w}(C_{m})=\mathrm{w}(C_{m}T_{1}^{\prime})$ is a
symplectic tableau hence $rC_{m}\leq lD_{m+1}$. So we can suppose
$\widetilde{e}_{i}C_{m}\neq0.$ Then $i\neq n$ and the column $C_{m}$ is
necessarily of type \textrm{(iv) }in (\ref{ei}). Write $C_{m}^{(i)}$ for the
admissible column such that \textrm{w}$(C_{m}^{(i)})=\widetilde{e}%
_{i}\mathrm{w}(C_{m})$.\ Then we have $rC_{m}=r(C_{m}^{(i)})$ and $rC_{m}\leq
lD_{m+1}$. Indeed by Formulas (\ref{TENS1}) and (\ref{TENS2}) $\widetilde
{e}_{i}^{r_{1}+1}\mathrm{w}(C_{m}T^{\prime})=\mathrm{w}(T^{\prime}%
)\widetilde{e}_{i}\mathrm{w}(C_{m})=$\textrm{w}$(C_{m}^{(i)}T^{\prime})$ and
$\widetilde{e}_{i}^{r_{1}+1}\mathrm{w}(C_{m}T^{\prime})$ is a symplectic tableau.

Now suppose that $C_{m}\neq D_{m}$ (hence $y\in C_{m}$ and $i\neq n$). We must
have $y=i+1\in C_{m},i\notin C_{m}$ and $\overline{i+1}\in C_{m}$ because
$C_{m}$ contains a letter of $\{\overline{i+1},i\}$. We have $r_{1}=p_{1}+1$
and $\widetilde{e}_{i}^{r_{1}}\mathrm{w}(C_{m}T^{\prime})=\mathrm{w}%
(T^{\prime})\,\widetilde{e}_{i}\mathrm{w}(C_{m})=\mathrm{w}(D_{m}T_{1}%
^{\prime})$ is a symplectic tableau. So $rD_{m}\leq lD_{m+1}$ and it suffices
to shows that $rC_{m-1}\leq lD_{m}$. The column $C_{m}$ is necessarily of type
\textrm{(ii) }or\textrm{ (v) }in \ref{ei} which implies that $lC_{m}%
=lC_{m}^{(i)}=lD_{m}$ where $C_{m}^{(i)}$ is the column with reading
$\widetilde{e}_{i}\mathrm{w}(C_{m})$. Then $rC_{m-1}\leq lC_{m}=lD_{m}$.
\end{proof}

\bigskip

\begin{example}
For $T=%
\begin{tabular}
[c]{|l|ll}\hline
$\mathtt{2}$ & $\mathtt{2}$ & \multicolumn{1}{|l|}{$\mathtt{3}$}\\\hline
$\mathtt{3}$ & $\mathtt{\bar{3}}$ & \multicolumn{1}{|l}{}\\\cline{1-1}%
\cline{1-2}%
$\mathtt{\bar{3}}$ &  & \\\cline{1-1}%
\end{tabular}
$ and $n=3$, we obtain successively\vspace{0.5cm}%

\begin{tabular}
[c]{|l|ll}\hline
$\mathtt{2}$ & $\mathtt{2}$ & \multicolumn{1}{|l|}{$\mathtt{2}$}\\\hline
$\mathtt{3}$ & $\mathtt{\bar{3}}$ & \multicolumn{1}{|l}{}\\\cline{1-1}%
\cline{1-2}%
$\mathtt{\bar{3}}$ &  & \\\cline{1-1}%
\end{tabular}
,
\begin{tabular}
[c]{|l|ll}\hline
$\mathtt{1}$ & $\mathtt{1}$ & \multicolumn{1}{|l|}{$\mathtt{1}$}\\\hline
$\mathtt{3}$ & $\mathtt{\bar{3}}$ & \multicolumn{1}{|l}{}\\\cline{1-1}%
\cline{1-2}%
$\mathtt{\bar{3}}$ &  & \\\cline{1-1}%
\end{tabular}
,
\begin{tabular}
[c]{|l|ll}\hline
$\mathtt{1}$ & $\mathtt{1}$ & \multicolumn{1}{|l|}{$\mathtt{1}$}\\\hline
$\mathtt{3}$ & $\mathtt{3}$ & \multicolumn{1}{|l}{}\\\cline{1-1}\cline{1-2}%
$\mathtt{\bar{3}}$ &  & \\\cline{1-1}%
\end{tabular}
,
\begin{tabular}
[c]{|l|ll}\hline
$\mathtt{1}$ & $\mathtt{1}$ & \multicolumn{1}{|l|}{$\mathtt{1}$}\\\hline
$\mathtt{2}$ & $\mathtt{2}$ & \multicolumn{1}{|l}{}\\\cline{1-1}\cline{1-2}%
$\mathtt{\bar{3}}$ &  & \\\cline{1-1}%
\end{tabular}
and
\begin{tabular}
[c]{|l|ll}\hline
$\mathtt{1}$ & $\mathtt{1}$ & \multicolumn{1}{|l|}{$\mathtt{1}$}\\\hline
$\mathtt{2}$ & $\mathtt{2}$ & \multicolumn{1}{|l}{}\\\cline{1-1}\cline{1-2}%
$\mathtt{3}$ &  & \\\cline{1-1}%
\end{tabular}
.%

\[
A(T)=f_{2}f_{1}^{(3)}f_{3}f_{2}^{(2)}f_{3}T_{\lambda}.
\]
\end{example}

\begin{proposition}
\label{prop_A(T)}The expansion of $A(T)$ on the basis $\{v_{\tau};\tau
\in\mathbf{T}(n,\lambda)\}$ of $W(\lambda)$ is of the form
\[
A(T)=\underset{\tau}{\sum}\alpha_{\tau,T}(q)v_{\tau}%
\]

where the coefficients $\alpha_{\tau,T}(q)$ satisfy:

\textrm{(i)}: $\alpha_{\tau,T}(q)\neq0$ only if $\tau$ and $T$ have the same weight,

\textrm{(ii)}: $\alpha_{\tau,T}(q)\in\mathbb{N}[q,q^{-1}]$ and $\alpha_{T,T}(q)=1,$

\textrm{(iii)}: $\alpha_{\tau,T}(q)\neq0$ only if $\tau\trianglelefteq T.$
\end{proposition}

\begin{proof}
\textrm{(i) }is a straightforward consequence of the definition of $A(T)$. By
Lemma \ref{Lem_F()}, we know that $\alpha_{\tau,T}(q)\in\mathbb{N}[q,q^{-1}].$
By induction the proposition will be proved if we show that \textrm{(ii)
}and\textrm{ (iii) }hold for $T$ as soon as they hold for $T_{1}$ with
$A(T)=f_{i_{1}}^{(r_{1})}A(T_{1})$ (the notations are those of Lemma
\ref{T1_Tab}). If the vector $v_{\tau}$ occurs in $A(T)$, the tabloid $\tau$
is obtained from a tabloid $\tau_{1}$ labelling a vector $v_{\tau_{1}}$
occurring in $A(T_{1})$ by changing $r_{1}$ letters $+$ of $\{\overline
{i+1},i\}$ into the corresponding letters $-$ of $\{\overline{i},i+1\}$.

It follows from the definition of $T_{1}$ that the tableau $T$ is obtained by
changing the $r_{1}$ rightmost letters $+$ of $T_{1}$ (i.e. the leftmost
letters $+$ of $\mathrm{w}(T_{1})$) into the corresponding letters $-$. Hence
$v_{T}$ appears in $f_{i_{1}}^{(r_{1})}v_{T_{1}}$ with a non zero coefficient.
Now suppose that there exists $\tau_{1}\neq T_{1}$ such that $v_{\tau_{1}}$
appears in $A(T_{1})$ and $v_{T}$ appears in $f_{i_{1}}^{(r_{1})}v_{\tau_{1}}%
$. Let $w$ be the factor of the words \textrm{w}$(\tau_{1})$ and
\textrm{w}$(T_{1})$ of maximal length such that there exist two words
$u,u^{\prime}$ and two letters $x\neq y$ satisfying:%
\begin{equation}
\mathrm{w}(\tau_{1})=wxu\text{ and }\mathrm{w}(T_{1})=wyu^{\prime} \label{dec}%
\end{equation}
We must have $x<y$ because $\tau_{1}\vartriangleleft T_{1}$. The letter $x$ is
necessarily modified when $T$ is obtained from $\tau_{1}$. Otherwise we have
$x\in$ \textrm{w}$(T)$. But \textrm{w}$(T)$ can also be computed from
\textrm{w}$(T_{1}).$ So the letter of \textrm{w}$(T)$ occurring at the same
place as the letter $y$ in \textrm{w}$(T_{1})$ is $\geq y$: it can not be $x$.
This implies that $x$ is letter $+$ and $y$ is its corresponding letter $-$.
Hence $w$ contains the $r_{1}$ letters $+$ changed to a $-$ when $T$ is
obtained from $T_{1}$. We derive a contradiction because in this case there
are $r_{1}+1$ letters $+$ changed into a $-$ when $T$ is obtained from
$\tau_{1}$. Hence $v_{T}$ can only appear in $f_{i_{1}}^{(r_{1})}v_{T_{1}}$.
With the notation of (\ref{coeff_N_sigma}) we will have $N(\mathrm{id}%
)=-\frac{r_{1}(r_{1}+1)}{2}$ because there is no letters $-$ to the left of
the letters $+$ modified in \textrm{w}$(T_{1})$ to obtain \textrm{w}$(T)$.
Then by (\ref{coeff_N_sigma}) the coefficient of $v_{T}$ in $f_{i_{1}}%
^{(r_{1})}v_{T_{1}}$ is equal to $1$ which proves $\mathrm{(ii)}$ for the
coefficient of $v_{T_{1}}$ in $A(T_{1})$ is $1$.

Consider $v_{\tau}$ appearing in $A(T)$ and suppose that the tabloid $\tau$ is
obtained from the tabloid $\tau_{1}$ such that $v_{\tau_{1}}$ appears in
$A(T_{1})$ by changing $r_{1}$ letters $+$ into their corresponding letters
$-$. Let $\tau^{\prime}$ be the tabloid obtained by changing in \textrm{w}%
$(\tau_{1})$ the $r_{1}$ leftmost letters $+$ (not immediately followed by
their corresponding letters $-$) by letters $-$. We are going to prove that
$\tau^{\prime}\trianglelefteq T$ which implies the proposition because
$\tau\trianglelefteq\tau^{\prime}$. If $\tau_{1}=T_{1}$ then $\tau^{\prime}%
=T$. So we can suppose $\tau_{1}\neq T_{1}$ and decompose the words
$\mathrm{w}(\tau_{1})$, $\mathrm{w}(T_{1})$ as in (\ref{dec}) with $x<y$. If
$\tau^{\prime}\trianglerighteq T,$ there is a letter $+$ in $w$ (that we write
$z_{+}$) which is changed into its corresponding letter $-$ when $\tau
^{\prime}$ is obtained from $\tau_{1}$ but is not modified when $T$ is
obtained from $T_{1}$. If we write $w=w_{1}z_{+}w_{2}$ where $w_{1}$ and
$w_{2}$ are words of $\mathcal{C}_{n}^{\ast}$, we have%
\[
\mathrm{w}(\tau_{1})=w_{1}z_{+}w_{2}xu\text{ and }\mathrm{w}(T_{1})=w_{1}%
z_{+}w_{2}yu^{\prime}.
\]
Then by definition of $T_{1}$, $w_{1}$ contains the $r_{1}$ letters $+$
changing to $-$ to obtain $T$. This contradicts the definition of
$\tau^{\prime}$. So \textrm{(iii) }is true.
\end{proof}

It follows from \textrm{(iii)} that the vectors $A(T)$ are linearly
independent in $V(\lambda)$. This implies that $\{A(T);T\in\mathbf{ST}%
(n,\lambda)\}$ is a $\mathbb{Q(}q\mathbb{)}$-basis of $V(\lambda)$. Indeed by
Theorem \ref{TH_KN}, we know that $\dim V(\lambda)=$\textrm{$card$%
}$(\mathbf{ST}(n,\lambda))$. As a consequence of (\ref{A(T)_monom}), we obtain
that $\overline{A(T)}=A(T)$. Note that, by definition of Marsh's algorithm,
the bases $\{A(T)\}$ and $\{G(T)\}$ coincide when $\lambda$ is a fundamental weight.

\subsection{From $\{A(T)\}$ to $\{G(T)\}$}

Let us write%
\[
G(T)=\underset{\tau\in\mathbf{T}(n,\lambda)}{\sum}d_{\tau,T}(q)\,v_{\tau}%
\]
We are going to describe a simple algorithm for computing the rectangular
matrix of coefficients
\[
D=[d_{\tau,T}(q)],\text{ \ \ \ \ }\tau\in\mathbf{T}(n,\lambda)\text{,
\ \ \ \ }T\in\mathbf{ST}(n,\lambda).
\]

\begin{lemma}
\label{G_on_tabloids}The coefficients $d_{\tau,T}(q)$ belong to $\mathbb{Q}%
[q]$. Moreover $d_{\tau,T}(0)=0$ if $\tau\neq T$ and $d_{T,T}(0)=1$.
\end{lemma}

\begin{proof}
Recall that $\{G(T)\}$ is a basis of $V_{\mathbb{Q}}(\lambda)=U_{\mathbb{Q}%
}^{-}v_{\lambda}$. This implies that the vectors of this basis are
$\mathbb{Q[}q,q^{-1}]$-linear combinations of vectors of the type considered
in Lemma \ref{Lem_F()}.\ In particular $d_{\tau,T}\in\mathbb{Q}[q,q^{-1}%
]$.\ By condition (\ref{cond_cong}), $d_{\tau,T}(q)$ must be regular at $q=0$
and
\[
d_{\tau,T}(q)\equiv\left\{
\begin{tabular}
[c]{l}%
$0$ $\operatorname{mod}q$ \ \ if $\tau\neq T$\\
$1$ $\operatorname{mod}q$ \ \ otherwise
\end{tabular}
\right.  .
\]
So $d_{\tau,T}(q)$ belong in fact to $\mathbb{Q[}q]$ and the Lemma is true.
\end{proof}

Let us write
\begin{equation}
G(T)=\underset{S\in ST(n,\lambda)}{\sum}\beta_{S,T}(q)\,A(S)
\label{G(T)_on_A(T)}%
\end{equation}
the expansion of the basis $\{G(T)\}$ on the basis $\{A(T)\}$. We have the
following lemma analogous to Lemma 4.3 of \cite{L-T}:

\begin{lemma}
The coefficients $\beta_{S,T}(q)$ of (\ref{G(T)_on_A(T)}) satisfy:

\textrm{(i)}: $\beta_{S,T}(q)=\beta_{S,T}(q^{-1}),$

\textrm{(ii)}: $\beta_{S,T}(q)=0$ unless $S\trianglelefteq T,$

\textrm{(iii)}: $\beta_{T,T}(q)=1$.
\end{lemma}

\begin{proof}
See proof of Lemma 4.3 in \cite{L-T}.
\end{proof}

Let $T_{\lambda}=T^{(1)}\vartriangleleft T^{(2)}\vartriangleleft\cdot
\cdot\cdot\vartriangleleft T^{(t)}$ be the sequence of tableaux of
$\mathbf{ST}(n,\lambda)$ ordered in increasing order. We have $G(T_{\lambda
})=A(T_{\lambda})$, i.e. $G(T^{(1)})=A(T^{(1)})$. By the previous lemma, the
transition matrix $M$ from $\{A(T)\}$ to $\{G(T)\}$ is upper unitriangular
once the two bases are ordered with $\trianglelefteq$. Since $\{G(T)\}$ is a
$\mathbb{Q}[q,q^{-1}]$ basis of $V_{\mathbb{Q}}(\lambda)$ and $A(T)\in
V_{\mathbb{Q}}(\lambda),$ the entries of $M$ are in $\mathbb{Q}[q,q^{-1}].$
Suppose by induction that we have computed the expansion on the basis
$\{v_{\tau};\tau\in\mathbf{T}(n,\lambda)\}$ of the vectors%
\[
G(T^{(1)}),...,G(T^{(i)})
\]
and that this expansion satisfies $d_{\tau,T^{(p)}}(q)=0$ if $\tau
\vartriangleright T^{(p)}$ for $p=1,...,i$. The inverse matrix $M^{-1}$ is
also upper unitriangular with entries in $\mathbb{Q}[q,q^{-1}]$. So we can
write:%
\begin{equation}
G(T^{(i+1)})=A(T^{(i+1)})-\gamma_{i}(q)G(T^{(i)})-\cdot\cdot\cdot-\gamma
_{1}(q)G(T^{(1)})\text{.} \label{G_on_A}%
\end{equation}
It follows from condition (\ref{cond_invo}) and Proposition \ref{prop_A(T)}
that $\gamma_{m}(q)=\gamma_{m}(q^{-1})$ for $m=1,...,i$. By Lemma
\ref{G_on_tabloids}, the coordinate $d_{T^{(i)},T^{(i+1)}}(q)$ of
$G(T^{(i+1)})$ on the vector $v_{T^{(i)}}$ belongs to $\mathbb{Q}[q]$,
$d_{T^{(i)},T^{(i+1)}}(0)=0$ and the coordinate $d_{T^{(i)},T^{(i)}}(q)$ of
$G(T^{(i)})$ on the vector $v_{T^{(i)}}$ is equal to $1$. Moreover
$v_{T^{(i)}}$ can only occur in $A(T^{(i+1)})-\gamma_{i}(q)G(T^{(i)}).$ If
\[
\alpha_{T^{(i)},T^{(i+1)}}(q)=\underset{j=-r}{\overset{s}{\sum}}a_{j}q^{j}%
\in\mathbb{N}[q,q^{-1}]
\]
then we will have
\[
\gamma_{i}(q)=\overset{0}{\underset{j=-r}{\sum}}a_{j}q^{j}+\underset
{j=1}{\overset{r}{\sum}}a_{-j}q^{j}\in\mathbb{N}[q,q^{-1}].
\]
Next if the coefficient of $v_{T^{(i-1)}}$ in $A(T^{(i+1)})-\gamma
_{i}(q)G(T^{(i)})$ is equal to
\[
\underset{j=-l}{\overset{k}{\sum}}b_{j}q^{j}%
\]
using similar arguments we obtain%
\[
\gamma_{i-1}(q)=\overset{0}{\underset{j=-l}{\sum}}b_{j}q^{j}+\overset
{l}{\underset{j=1}{\sum}}b_{-j}q^{j},
\]
and so on. So we have computed the expansion of $G(T^{(i+1)})$ on the basis
$\{v_{\tau}\}$ and this expansion satisfies $d_{\tau,T^{(i+1)}}(q)=0$ if
$\tau\vartriangleright T^{(i+1)}$. Finally notice that $\gamma_{s}%
(q)\in\mathbb{Z}[q,q^{-1}]$ for all $s$ by Proposition \ref{prop_A(T)}.

\noindent By construction of $A(T)$, it is possible to write the basis
$\{A(T)\}$ in terms of the basis $\{v_{\tau}\}.\;$Using the above, it is then
possible to determine the $\gamma_{i}(q)$ and thus write the $G(T)$ in terms
of the basis $\{v_{\tau}\}.$ We have proved that:

\begin{theorem}
Let $T\in\mathbf{ST}(n,\lambda).\;$Then $G(T)=\sum d_{\tau,T}(q)v_{\tau}$
where the coefficients $d_{\tau,T}(q)$ satisfy:

\textrm{(i)}: $d_{\tau,T}(q)\in\mathbb{Z}[q],$

\textrm{(ii)}: $d_{T,T}(q)=1$ and $d_{\tau,T}(0)=0$ for $\tau\neq T,$

\textrm{(iii)}: $d_{\tau,T}(q)\neq0$ only if $\tau$ and $T$ have the same
weight, and $\tau\trianglelefteq T$.
\end{theorem}

\section{Examples}

All the vectors occurring in our calculations are weight vectors.\ So we can
use our algorithm to compute the canonical basis of a single weight space. We
give below the matrix obtained for the $12$-dimension weight space of the
$U_{q}(sp_{6})$-module $V(4,3,2)$ (i.e. $\lambda=\Lambda_{1}+\Lambda
_{2}+2\Lambda_{3}$) corresponding to the weight $\mu=(0,3,0).$ Its columns and
rows are respectively labelled by the symplectic tableaux and by the tabloids
of weight $\mu$ ordered from left to right and top to bottom in decreasing
order for $\trianglelefteq$.\ Those tabloids which are symplectic tableaux
have been written in bold style.

\bigskip

{\scriptsize \hskip-16mm $
%
%
%
%
%
%
%
%
%
%
%
%
%
%
%
%
%
%
%
%
%
%
%
%
%
%
%
%
%
%
%
%
%
%
%
%
%
%
%
%
%
%
%
%
%
%
%
%
%
%
%
%
%
%
%
%
%
%
%
%
%
%
%
\begin{array}
[c]{ccccccccccccc}%
. &
\begin{array}
[c]{c}%
{11\bar{2}\bar{1}}\text{ }\\
{3\bar{3}\bar{1}{\ \ }}\\
{\bar{2}\bar{2}{\ }{\ }{\ }}%
\end{array}
{\ } &
\begin{array}
[c]{c}%
{13\bar{3}\bar{1}}\text{ }\\
{3\bar{3}\bar{2}{\ \ }}\\
{\bar{2}\bar{2}{\ }{\ }{\ }}%
\end{array}
{\ } &
\begin{array}
[c]{c}%
{133\bar{1}}\text{ }\\
{\bar{3}\bar{3}\bar{2}{\ \ }}\\
{\bar{2}\bar{2}{\ }{\ }{\ }}%
\end{array}
{\ } &
\begin{array}
[c]{c}%
{13\bar{2}\bar{2}}\text{ }\\
{3\bar{3}\bar{1}{\ \ }}\\
{\bar{3}\bar{2}{\ }{\ }{\ }}%
\end{array}
{\ } &
\begin{array}
[c]{c}%
{12\bar{2}\bar{2}}\text{ }\\
{3\bar{3}\bar{1}{\ \ }}\\
{\bar{2}\bar{2}{\ }{\ }{\ }}%
\end{array}
{\ } &
\begin{array}
[c]{c}%
{13\bar{3}\bar{2}}\text{ }\\
{3\bar{3}\bar{1}{\ \ }}\\
{\bar{2}\bar{2}{\ }{\ }{\ }}%
\end{array}
{\ } &
\begin{array}
[c]{c}%
{13\bar{3}\bar{2}}\text{ }\\
{3\bar{2}\bar{2}{\ \ }}\\
{\bar{3}\bar{1}{\ }{\ }{\ }}%
\end{array}
{\ } &
\begin{array}
[c]{c}%
{13\bar{3}\bar{2}}\text{ }\\
{3\bar{3}\bar{2}{\ \ }}\\
{\bar{2}\bar{1}{\ }{\ }{\ }}%
\end{array}
{\ } &
\begin{array}
[c]{c}%
{133\bar{2}}\text{ }\\
{\bar{3}\bar{3}\bar{1}{\ \ }}\\
{\bar{2}\bar{2}{\ }{\ }{\ }}%
\end{array}
{\ } &
\begin{array}
[c]{c}%
{133\bar{2}}\text{ }\\
{\bar{3}\bar{3}\bar{2}{\ \ }}\\
{\bar{2}\bar{1}{\ }{\ }{\ }}%
\end{array}
{\ } &
\begin{array}
[c]{c}%
{13\bar{3}\bar{3}}\text{ }\\
{3\bar{2}\bar{2}{\ \ }}\\
{\bar{2}\bar{1}{\ }{\ \ }}%
\end{array}
{\ } &
\begin{array}
[c]{c}%
{133\bar{3}}\text{ }\\
{\bar{3}\bar{2}\bar{2}{\ \ }}\\
{\bar{2}\bar{1}{\ }{\ }{\ }}%
\end{array}
{\ }\\
{\
\begin{array}
[c]{c}%
\mathbf{11\bar{2}\bar{1}}\text{ }\\
\mathbf{3\bar{3}\bar{1}\ \ }\\
\mathbf{\bar{2}\bar{2}\ \ \ }%
\end{array}
} & 1 & . & . & . & . & . & . & . & . & . & . & .\\
{\ }%
\begin{array}
[c]{c}%
11\bar{2}\bar{1}\\
\bar{3}3\bar{1}\ \ \\
\bar{2}\bar{2}\ \ \
\end{array}
& q^{2} & . & . & . & . & . & . & . & . & . & . & .\\
{\ }%
\begin{array}
[c]{c}%
\mathbf{13\bar{3}\bar{1}}\text{ }\\
\mathbf{3\bar{3}\bar{2}\ \ }\\
\mathbf{\bar{2}\bar{2}\ \ \ }%
\end{array}
& q & 1 & . & . & . & . & . & . & . & . & . & .\\
{\ }%
\begin{array}
[c]{c}%
11\bar{3}\bar{1}{{\ }}\\
3\bar{2}\bar{2}\ \ \\
\bar{2}\bar{1}\ \ \
\end{array}
& q^{2} & q & . & . & . & . & . & . & . & . & . & .\\
{\ }%
\begin{array}
[c]{c}%
{31\bar{3}\bar{1}{\ }}\\
{\bar{3}3\bar{2}{\ \ }}\\
{\bar{2}\bar{2}{\ }{\ }{\ }}%
\end{array}
& q^{3} & q^{2} & . & . & . & . & . & . & . & . & . & .\\
{\ }%
\begin{array}
[c]{c}%
{11\bar{3}\bar{1}{\ }}\\
{\bar{2}3\bar{2}{\ \ }}\\
{\bar{1}\bar{2}{\ }{\ }{\ }}%
\end{array}
& q^{4} & q^{3} & . & . & . & . & . & . & . & . & . & .\\
{\ }%
\begin{array}
[c]{c}%
\mathbf{133\bar{1}}{{\ }}\\
\mathbf{\bar{3}\bar{3}\bar{2}\ \ }\\
\mathbf{\bar{2}\bar{2}\ \ \ }%
\end{array}
& q^{3} & q^{2} & 1 & . & . & . & . & . & . & . & . & .\\
{\ }%
\begin{array}
[c]{c}%
{113\bar{1}{\ }}\\
{\bar{3}\bar{2}\bar{2}{\ \ }}\\
{\bar{2}1{\ }{\ }{\ }}%
\end{array}
& q^{4} & q^{3} & q & . & . & . & . & . & . & . & . & .\\
{\ }%
\begin{array}
[c]{c}%
{313\bar{1}{\ }}\\
{\bar{3}\bar{3}\bar{2}{\ \ }}\\
{\bar{2}\bar{2}{\ }{\ }{\ }}%
\end{array}
& q^{5} & q^{4} & q^{2} & {.} & . & . & . & . & . & . & . & .\\
{\ }%
\begin{array}
[c]{c}%
{113\bar{1}{\ }}\\
{\bar{2}\bar{3}\bar{2}{\ \ }}\\
{\bar{1}\bar{2}{\ }{\ }{\ }}%
\end{array}
& q^{6} & q^{5} & q^{3} & . & . & . & . & . & . & . & . & .\\
{\ }%
\begin{array}
[c]{c}%
{1\bar{3}1\bar{1}{\ }}\\
{3\bar{2}\bar{2}{\ \ }}\\
{\bar{2}\bar{1}{\ }{\ }{\ }}%
\end{array}
& {.} & q^{2} & {.} & {.} & . & . & . & . & . & . & . & .\\
{\ }%
\begin{array}
[c]{c}%
{1311{\ }}\\
{\bar{3}\bar{2}\bar{2}{\ \ }}\\
{\bar{2}\bar{1}{\ }{\ }{\ }}%
\end{array}
& {.} & q^{4} & q^{2} & . & {.} & . & . & . & . & . & . & .\\
{\ }%
\begin{array}
[c]{c}%
{331\bar{1}{\ }}\\
{\bar{3}\bar{3}\bar{2}{\ \ }}\\
{\bar{2}\bar{2}{\ }{\ }{\ }}%
\end{array}
& . & . & q^{3} & . & . & . & . & . & . & . & . & .\\
{\ }%
\begin{array}
[c]{c}%
{1311{\ }}\\
{\bar{2}\bar{3}\bar{2}{\ \ }}\\
{\bar{1}\bar{2}{\ }{\ }{\ }}%
\end{array}
& . & . & {q}^{4} & . & . & . & . & . & . & . & . & .\\
{\ }%
\begin{array}
[c]{c}%
{3111{\ }}\\
{\bar{3}\bar{2}\bar{2}{\ \ }}\\
{\bar{2}\bar{1}{\ }{\ }{\ }}%
\end{array}
& . & . & q^{4} & {.} & . & . & . & . & . & . & . & .\\
{\ }%
\begin{array}
[c]{c}%
{111\bar{1}{\ }}\\
{\bar{2}\bar{2}\bar{2}{\ \ }}\\
{\bar{1}\bar{1}{\ }{\ }{\ }}%
\end{array}
& . & . & {q}^{5} & . & . & . & . & . & . & . & . & .\\
{\ }%
\begin{array}
[c]{c}%
{3111{\ }}\\
{\bar{2}\bar{3}\bar{2}{\ \ }}\\
{\bar{1}\bar{2}{\ }{\ }{\ }}%
\end{array}
& {.} & q^{4} & q^{6} & {.} & . & . & . & . & . & . & . & .\\
{\ }%
\begin{array}
[c]{c}%
{\bar{3}11\bar{1}{\ }}\\
{\bar{2}3\bar{2}{\ \ }}\\
{\bar{1}\bar{2}{\ }{\ }{\ }}%
\end{array}
& . & q^{6} & . & {.} & {.} & . & . & . & . & . & . & .\\
{\ }%
\begin{array}
[c]{c}%
\mathbf{13\bar{2}\bar{2}}{{\ }}\\
\mathbf{3\bar{3}\bar{1}\ \ }\\
\mathbf{\bar{3}\bar{2}\ \ \ }%
\end{array}
& {.} & . & {.} & {1} & . & . & . & . & . & . & . & .\\
{\ }%
\begin{array}
[c]{c}%
{13\bar{2}\bar{2}{\ }}\\
{2\bar{3}\bar{1}{\ \ }}\\
{\bar{2}\bar{2}{\ }{\ }{\ }}%
\end{array}
& {.} & {.} & {.} & q & {.} & . & . & . & . & . & . & .\\
{\ }%
\begin{array}
[c]{c}%
\mathbf{12\bar{2}\bar{2}}{{\ }}\\
\mathbf{3\bar{3}\bar{1}\ \ }\\
\mathbf{\bar{2}\bar{2}\ \ \ }%
\end{array}
& q & . & . & q^{2} & 1 & . & . & . & . & . & . & .\\
{\ }%
\begin{array}
[c]{c}%
{12\bar{2}\bar{2}{\ }}\\
{\bar{3}3\bar{1}{\ \ }}\\
{\bar{2}\bar{2}{\ }{\ }{\ }}%
\end{array}
& q^{3} & . & . & . & q^{2} & . & . & . & . & . & . & .\\
{\ }%
\begin{array}
[c]{c}%
{11\bar{2}\bar{2}{\ }}\\
{3\bar{2}1{\ \ }}\\
{\bar{3}\bar{1}{\ }{\ }{\ }}%
\end{array}
& . & . & . & q & . & . & . & . & . & . & . & .\\
{\ }%
\begin{array}
[c]{c}%
{11\bar{2}\bar{2}{\ }}\\
{2\bar{2}\bar{1}{\ \ }}\\
{\bar{2}\bar{1}{\ }{\ }{\ }}%
\end{array}
& . & . & . & q^{2} & . & . & . & . & . & . & . & .\\
{\ }%
\begin{array}
[c]{c}%
{11\bar{2}\bar{2}{\ }}\\
{3\bar{3}\bar{1}{\ \ }}\\
{\bar{2}\bar{1}{\ }{\ }{\ }}%
\end{array}
& q^{2} & . & . & q^{3} & q & . & . & . & . & . & . & .
\end{array}
%
%
%
%
%
%
%
%
%
%
%
%
%
%
%
%
%
%
%
%
%
%
%
%
%
%
%
%
%
%
%
%
%
%
%
%
%
%
%
%
%
%
%
%
%
%
%
%
%
%
%
%
%
%
%
%
%
%
%
%
%
%
%
$ }

{\scriptsize \hskip-16mm $
%
%
%
%
%
%
%
%
%
%
%
%
%
%
%
%
%
%
%
%
%
%
%
%
%
%
%
%
%
%
%
%
%
%
%
%
%
%
%
%
%
%
%
%
%
%
%
%
%
%
%
%
%
%
%
%
%
%
%
%
%
%
%
\begin{array}
[c]{ccccccccccccc}%
&
\begin{array}
[c]{c}%
{11\bar{2}\bar{1}}\text{ }\\
{3\bar{3}\bar{1}{\ \ }}\\
{\bar{2}\bar{2}{\ }{\ }{\ }}%
\end{array}
&
\begin{array}
[c]{c}%
{13\bar{3}\bar{1}}\text{ }\\
{3\bar{3}\bar{2}{\ \ }}\\
{\bar{2}\bar{2}{\ }{\ }{\ }}%
\end{array}
&
\begin{array}
[c]{c}%
{133\bar{1}}\text{ }\\
{\bar{3}\bar{3}\bar{2}{\ \ }}\\
{\bar{2}\bar{2}{\ }{\ }{\ }}%
\end{array}
&
\begin{array}
[c]{c}%
{13\bar{2}\bar{2}}\text{ }\\
{3\bar{3}\bar{1}{\ \ }}\\
{\bar{3}\bar{2}{\ }{\ }{\ }}%
\end{array}
&
\begin{array}
[c]{c}%
{12\bar{2}\bar{2}}\text{ }\\
{3\bar{3}\bar{1}{\ \ }}\\
{\bar{2}\bar{2}{\ }{\ }{\ }}%
\end{array}
&
\begin{array}
[c]{c}%
{13\bar{3}\bar{2}}\text{ }\\
{3\bar{3}\bar{1}{\ \ }}\\
{\bar{2}\bar{2}{\ }{\ }{\ }}%
\end{array}
&
\begin{array}
[c]{c}%
{13\bar{3}\bar{2}}\text{ }\\
{3\bar{2}\bar{2}{\ \ }}\\
{\bar{3}\bar{1}{\ }{\ }{\ }}%
\end{array}
&
\begin{array}
[c]{c}%
{13\bar{3}\bar{2}}\text{ }\\
{3\bar{3}\bar{2}{\ \ }}\\
{\bar{2}\bar{1}{\ }{\ }{\ }}%
\end{array}
&
\begin{array}
[c]{c}%
{133\bar{2}}\text{ }\\
{\bar{3}\bar{3}\bar{1}{\ \ }}\\
{\bar{2}\bar{2}{\ }{\ }{\ }}%
\end{array}
&
\begin{array}
[c]{c}%
{133\bar{2}}\text{ }\\
{\bar{3}\bar{3}\bar{2}{\ \ }}\\
{\bar{2}\bar{1}{\ }{\ }{\ }}%
\end{array}
&
\begin{array}
[c]{c}%
{13\bar{3}\bar{3}}\text{ }\\
{3\bar{2}\bar{2}{\ \ }}\\
{\bar{2}\bar{1}{\ }{\ \ }}%
\end{array}
&
\begin{array}
[c]{c}%
{133\bar{3}}\text{ }\\
{\bar{3}\bar{2}\bar{2}{\ \ }}\\
{\bar{2}\bar{1}{\ }{\ }{\ }}%
\end{array}
\\
{\
\begin{array}
[c]{c}%
21\bar{2}\bar{2}{{\ }}\\
3\bar{3}\bar{1}\ \ \\
\bar{2}\bar{2}\ \ \
\end{array}
} & q^{3} & . & . & . & q^{2} & . & . & . & . & . & . & .\\
{\
\begin{array}
[c]{c}%
11\bar{2}\bar{2}{{\ }}\\
3\bar{3}\bar{1}\ \ \\
\bar{1}\bar{2}\ \ \
\end{array}
} & q^{4} & . & . & . & q^{3} & . & . & . & . & . & . & .\\
{\
\begin{array}
[c]{c}%
21\bar{2}\bar{2}{{\ }}\\
\bar{3}3\bar{1}\ \ \\
\bar{2}\bar{1}\ \ \
\end{array}
} & q^{4} & . & . & . & q^{3} & . & . & . & . & . & . & .\\
{\ }%
\begin{array}
[c]{c}%
21\bar{2}\bar{2}{{\ }}\\
\bar{3}3\bar{1}\ \ \\
\bar{2}\bar{2}\ \ \
\end{array}
& q^{5} & . & . & q^{2} & .q^{4} & . & . & . & . & . & . & .\\
{\ }%
\begin{array}
[c]{c}%
11\bar{2}\bar{2}{{\ }}\\
\bar{3}3\bar{1}\ \ \\
\bar{1}\bar{2}\ \ \
\end{array}
& q^{6} & . & . & q^{3} & q^{5} & . & . & . & . & . & . & .\\
{\ }%
\begin{array}
[c]{c}%
31\bar{2}\bar{2}{{\ }}\\
\bar{3}3\bar{1}\ \ \\
\bar{2}\bar{3}\ \ \
\end{array}
& . & {.} & . & q^{3} & . & . & . & . & . & . & . & .\\
{\ }%
\begin{array}
[c]{c}%
11\bar{2}\bar{2}{{\ }}\\
\bar{2}3\bar{1}\ \ \\
\bar{1}\bar{3}\ \ \
\end{array}
& . & . & . & q^{4} & . & . & . & . & . & . & . & .\\
{\ }%
\begin{array}
[c]{c}%
{31\bar{2}\bar{2}{\ }}\\
{\bar{3}2\bar{1}{\ \ }}\\
{\bar{2}\bar{2}{\ }{\ }{\ }}%
\end{array}
& . & {.} & {.} & q^{4} & . & . & . & . & . & . & . & .\\
{\ }%
\begin{array}
[c]{c}%
{11\bar{2}\bar{2}{\ }}\\
{\bar{2}2\bar{1}{\ \ }}\\
{\bar{1}\bar{2}{\ }{\ }{\ }}%
\end{array}
& {.} & {.} & . & {q}^{5} & . & . & . & . & . & . & . & .\\
{\ }%
\begin{array}
[c]{c}%
\mathbf{13\bar{3}\bar{2}}{{\ }}\\
\mathbf{3\bar{3}\bar{1}\ \ }\\
\mathbf{\bar{2}\bar{2}\ \ \ }%
\end{array}
& q^{2} & q & {.} & q^{3} & q & 1 & . & . & . & . & . & .\\
{\ }%
\begin{array}
[c]{c}%
{11\bar{3}\bar{2}{\ }}\\
{3\bar{2}1{\ \ }}\\
{\bar{2}\bar{1}{\ }{\ }{\ }}%
\end{array}
& q^{3} & q^{2} & {.} & {q}^{4} & q^{2} & q & . & . & . & . & . & .\\
{\ }%
\begin{array}
[c]{c}%
{31\bar{3}\bar{2}{\ }}\\
{\bar{3}3\bar{1}{\ \ }}\\
{\bar{2}\bar{2}{\ }{\ }{\ }}%
\end{array}
& {q}^{4} & q^{3} & . & q^{5} & q^{3} & q^{2} & . & . & . & . & . & .\\
{\ }%
\begin{array}
[c]{c}%
{11\bar{3}\bar{2}{\ }}\\
{\bar{2}3\bar{1}{\ \ }}\\
{\bar{1}\bar{2}{\ }{\ }{\ }}%
\end{array}
& q^{5} & q^{4} & . & q^{6} & q^{4} & q^{3} & . & . & . & . & . & .\\
{\ }%
\begin{array}
[c]{c}%
\mathbf{13\bar{3}\bar{2}}{{\ }}\\
\mathbf{3\bar{2}\bar{2}\ \ }\\
\mathbf{\bar{3}\bar{1}\ \ \ }%
\end{array}
& . & . & {.} & q^{2} & . & . & 1 & . & . & . & . & .\\
{\ }%
\begin{array}
[c]{c}%
{13\bar{3}\bar{2}{\ }}\\
{2\bar{2}\bar{2}{\ \ }}\\
{\bar{2}\bar{1}{\ }{\ }{\ }}%
\end{array}
& . & . & . & q^{3} & . & . & q & . & . & . & . & .\\
{\ }%
\begin{array}
[c]{c}%
\mathbf{13\bar{3}\bar{2}}{{\ }}\\
\mathbf{3\bar{3}\bar{2}\ \ }\\
\mathbf{\bar{2}\bar{1}\ \ \ }%
\end{array}
& q^{3} & q^{2} & {.} & q^{4} & q^{2} & q & q^{2} & 1 & . & . & . & .\\
{\ }%
\begin{array}
[c]{c}%
{23\bar{3}\bar{2}{\ }}\\
{3\bar{3}\bar{2}{\ \ }}\\
{\bar{2}\bar{2}{\ }{\ }{\ }}%
\end{array}
& {q}^{4} & q^{3} & {.} & {.} & q^{3} & q^{2} & . & q & . & . & . & .\\
{\ }%
\begin{array}
[c]{c}%
{13\bar{3}\bar{2}{\ }}\\
{3\bar{3}\bar{2}{\ \ }}\\
{\bar{1}\bar{2}{\ }{\ }{\ }}%
\end{array}
& q^{5} & q^{4} & . & {.} & q^{4} & q^{3} & . & q^{2} & . & . & . & .\\
{\ }%
\begin{array}
[c]{c}%
{12\bar{3}\bar{2}{\ }}\\
{3\bar{2}\bar{2}{\ \ }}\\
{\bar{2}\bar{1}{\ }{\ }{\ }}%
\end{array}
& q^{4} & q^{3} & {.} & q^{5} & q^{3} & q^{2} & q^{3} & q & . & . & . & .\\
{\ }%
\begin{array}
[c]{c}%
{32\bar{3}\bar{2}{\ }}\\
{\bar{3}3\bar{2}{\ \ }}\\
{\bar{2}\bar{2}{\ }{\ }{\ }}%
\end{array}
& {q}^{5} & q^{4} & {.} & {.} & q^{4} & q^{3} & . & q^{2} & . & . & . & .\\
{\ }%
\begin{array}
[c]{c}%
{12\bar{3}\bar{2}{\ }}\\
{\bar{2}3\bar{2}{\ \ }}\\
{\bar{1}\bar{2}{\ }{\ }{\ }}%
\end{array}
& q^{6} & q^{5} & . & . & q^{5} & q^{4} & . & q^{3} & . & . & . & .\\
{\ }%
\begin{array}
[c]{c}%
{21\bar{3}\bar{2}{\ }}\\
{3\bar{2}\bar{2}{\ \ }}\\
{\bar{2}\bar{1}{\ }{\ }{\ }}%
\end{array}
& q^{5} & q^{4} & . & . & q^{4} & q^{3} & . & q^{2} & . & . & . & .\\
{\ }%
\begin{array}
[c]{c}%
{11\bar{3}\bar{2}{\ }}\\
{3\bar{2}\bar{2}{\ \ }}\\
{\bar{1}\bar{1}{\ }{\ }{\ }}%
\end{array}
& q^{6} & q^{5} & . & . & q^{5} & q^{4} & . & q^{3} & . & . & . & .\\
{\ }%
\begin{array}
[c]{c}%
{31\bar{3}\bar{2}{\ }}\\
{\bar{3}3\bar{2}{\ \ }}\\
{\bar{2}\bar{1}{\ }{\ }{\ }}%
\end{array}
& q^{6} & q^{5} & . & . & q^{5} & q^{4} & . & q^{3} & . & . & . & .\\
{\ }%
\begin{array}
[c]{c}%
{11\bar{3}\bar{2}{\ }}\\
{\bar{2}3\bar{2}{\ \ }}\\
{\bar{1}\bar{1}{\ }{\ }{\ }}%
\end{array}
& q^{7} & q^{6} & . & . & q^{6} & q^{5} & . & q^{4} & . & . & . & .
\end{array}
%
%
%
%
%
%
%
%
%
%
%
%
%
%
%
%
%
%
%
%
%
%
%
%
%
%
%
%
%
%
%
%
%
%
%
%
%
%
%
%
%
%
%
%
%
%
%
%
%
%
%
%
%
%
%
%
%
%
%
%
%
%
%
$}

\bigskip

{\scriptsize \hskip-16mm $
%
%
%
%
%
%
%
%
%
%
%
%
%
%
%
%
%
%
%
%
%
%
%
%
%
%
%
%
%
%
%
%
%
%
%
%
%
%
%
%
%
%
%
%
%
%
%
%
%
%
%
%
%
%
%
%
%
%
%
%
%
%
%
%
%
%
%
%
%
\begin{array}
[c]{ccccccccccccc}%
. & {\ }%
\begin{array}
[c]{c}%
{11\bar{2}\bar{1}}\text{ }\\
{3\bar{3}\bar{1}{\ \ }}\\
{\bar{2}\bar{2}{\ }{\ }{\ }}%
\end{array}
&
\begin{array}
[c]{c}%
{13\bar{3}\bar{1}}\text{ }\\
{3\bar{3}\bar{2}{\ \ }}\\
{\bar{2}\bar{2}{\ }{\ }{\ }}%
\end{array}
{\ } &
\begin{array}
[c]{c}%
{133\bar{1}}\text{ }\\
{\bar{3}\bar{3}\bar{2}{\ \ }}\\
{\bar{2}\bar{2}{\ }{\ }{\ }}%
\end{array}
{\ } &
\begin{array}
[c]{c}%
{13\bar{2}\bar{2}}\text{ }\\
{3\bar{3}\bar{1}{\ \ }}\\
{\bar{3}\bar{2}{\ }{\ }{\ }}%
\end{array}
{\ } & {\ }%
\begin{array}
[c]{c}%
{12\bar{2}\bar{2}}\text{ }\\
{3\bar{3}\bar{1}{\ \ }}\\
{\bar{2}\bar{2}{\ }{\ }{\ }}%
\end{array}
& {\ }%
\begin{array}
[c]{c}%
{13\bar{3}\bar{2}}\text{ }\\
{3\bar{3}\bar{1}{\ \ }}\\
{\bar{2}\bar{2}{\ }{\ }{\ }}%
\end{array}
&
\begin{array}
[c]{c}%
{13\bar{3}\bar{2}}\text{ }\\
{3\bar{2}\bar{2}{\ \ }}\\
{\bar{3}\bar{1}{\ }{\ }{\ }}%
\end{array}
{\ } & {\ }%
\begin{array}
[c]{c}%
{13\bar{3}\bar{2}}\text{ }\\
{3\bar{3}\bar{2}{\ \ }}\\
{\bar{2}\bar{1}{\ }{\ }{\ }}%
\end{array}
& {\ }%
\begin{array}
[c]{c}%
{133\bar{2}}\text{ }\\
{\bar{3}\bar{3}\bar{1}{\ \ }}\\
{\bar{2}\bar{2}{\ }{\ }{\ }}%
\end{array}
&
\begin{array}
[c]{c}%
{133\bar{2}}\text{ }\\
{\bar{3}\bar{3}\bar{2}{\ \ }}\\
{\bar{2}\bar{1}{\ }{\ }{\ }}%
\end{array}
{\ } & {\ }%
\begin{array}
[c]{c}%
{13\bar{3}\bar{3}}\text{ }\\
{3\bar{2}\bar{2}{\ \ }}\\
{\bar{2}\bar{1}{\ }{\ }{\ }}%
\end{array}
&
\begin{array}
[c]{c}%
{133\bar{3}}\text{ }\\
{\bar{3}\bar{2}\bar{2}{\ \ }}\\
{\bar{2}\bar{1}{\ }{\ }{\ }}%
\end{array}
{\ }\\
{\ }%
\begin{array}
[c]{c}%
{31\bar{3}\bar{2}{\ }}\\
{\bar{3}3\bar{2}{\ \ }}\\
{\bar{1}\bar{2}{\ }{\ }{\ }}%
\end{array}
& q^{7} & q^{6} & . & q^{4} & q^{6} & q^{5} & q^{2} & q^{4} & . & . & . & .\\
{\ }%
\begin{array}
[c]{c}%
{21\bar{3}\bar{2}{\ }}\\
{\bar{2}3\bar{2}{\ \ }}\\
{\bar{1}\bar{2}{\ }{\ }{\ }}%
\end{array}
& q^{8} & q^{7} & . & q^{5} & q^{7} & q^{6} & q^{3} & q^{5} & . & . & . & .\\
{\ }%
\begin{array}
[c]{c}%
{31\bar{3}\bar{2}{\ }}\\
{\bar{2}3\bar{2}{\ \ }}\\
{\bar{1}\bar{3}{\ }{\ }{\ }}%
\end{array}
& . & . & . & q^{6} & . & . & q^{4} & . & . & . & . & .\\
{\ }%
\begin{array}
[c]{c}%
{31\bar{3}\bar{2}{\ }}\\
{\bar{2}2\bar{2}{\ \ }}\\
{\bar{1}\bar{2}{\ }{\ }{\ }}%
\end{array}
& . & . & . & q^{7} & . & . & q^{5} & . & . & . & . & .\\
{\ }%
\begin{array}
[c]{c}%
\mathbf{133\bar{2}}{{\ }}\\
\mathbf{\bar{3}\bar{3}\bar{1}\ \ }\\
\mathbf{\bar{2}\bar{2}\ \ \ }%
\end{array}
& q^{4} & q^{3} & q & . & q^{3} & q^{2} & . & . & 1 & . & . & .\\
{\ }%
\begin{array}
[c]{c}%
{113\bar{2}{\ }}\\
{\bar{3}\bar{2}\bar{1}{\ \ }}\\
{\bar{2}\bar{1}{\ }{\ }{\ }}%
\end{array}
& q^{5} & q^{4} & q^{2} & . & q^{4} & q^{3} & . & . & q & . & . & .\\
{\ }%
\begin{array}
[c]{c}%
{313\bar{2}{\ }}\\
{\bar{3}\bar{3}\bar{1}{\ \ }}\\
{\bar{2}\bar{2}{\ }{\ }{\ }}%
\end{array}
& q^{6} & q^{5} & q^{2} & . & q^{5} & q^{4} & . & . & q^{2} & . & . & .\\
{\ }%
\begin{array}
[c]{c}%
{113\bar{2}{\ }}\\
{\bar{2}\bar{3}\bar{1}{\ \ }}\\
{\bar{1}\bar{2}{\ }{\ }{\ }}%
\end{array}
& q^{7} & q^{6} & q^{4} & . & q^{6} & q^{5} & . & . & q^{3} & . & . & .\\
{\ }%
\begin{array}
[c]{c}%
{1\bar{3}3\bar{2}{\ }}\\
{3\bar{2}\bar{2}{\ \ }}\\
{\bar{3}\bar{1}{\ }{\ }{\ }}%
\end{array}
& {.} & {.} & . & {.} & . & . & q^{2} & . & . & . & . & .\\
{\ }%
\begin{array}
[c]{c}%
{1\bar{3}3\bar{2}{\ }}\\
{2\bar{2}\bar{2}{\ \ }}\\
{\bar{2}\bar{1}{\ }{\ }{\ }}%
\end{array}
& . & . & {.} & . & . & . & q^{3} & . & . & . & . & .\\
{\ }%
\begin{array}
[c]{c}%
\mathbf{133\bar{2}}{{\ }}\\
\mathbf{\bar{3}\bar{3}\bar{2}\ \ }\\
\mathbf{\bar{2}\bar{1}\ \ \ }%
\end{array}
& q^{5} & q^{4} & q^{2} & q^{2} & q^{4} & q^{3} & q^{4} & q^{2} & q & 1 & . &
.\\
{\ }%
\begin{array}
[c]{c}%
{233\bar{2}{\ }}\\
{\bar{3}\bar{3}\bar{2}{\ \ }}\\
{\bar{2}\bar{2}{\ }{\ }{\ }}%
\end{array}
& q^{6} & q^{5} & q^{3} & q^{3} & q^{5} & q^{4} & . & q^{3} & q^{2} & q & . &
.\\
{\ }%
\begin{array}
[c]{c}%
{133\bar{2}{\ }}\\
{\bar{3}\bar{3}\bar{2}{\ \ }}\\
{\bar{1}\bar{2}{\ }{\ }{\ }}%
\end{array}
& q^{7} & q^{6} & q^{4} & q^{4} & q^{6} & q^{5} & . & q^{4} & q^{3} & q^{2} &
. & .\\
{\ }%
\begin{array}
[c]{c}%
{123\bar{2}{\ }}\\
{\bar{3}\bar{2}\bar{2}{\ \ }}\\
{\bar{2}\bar{1}{\ }{\ }{\ }}%
\end{array}
& q^{6} & q^{5} & q^{3} & q^{3} & q^{5} & q^{4} & q^{5} & q^{3} & q^{2} & q &
. & .\\
{\ }%
\begin{array}
[c]{c}%
{323\bar{2}{\ }}\\
{\bar{3}\bar{3}\bar{2}{\ \ }}\\
{\bar{2}\bar{2}{\ }{\ }{\ }}%
\end{array}
& q^{7} & q^{6} & q^{4} & q^{4} & q^{6} & q^{5} & . & q^{4} & q^{3} & q^{2} &
. & .\\
{\ }%
\begin{array}
[c]{c}%
{123\bar{2}{\ }}\\
{\bar{2}\bar{3}\bar{2}{\ \ }}\\
{\bar{1}\bar{2}{\ }{\ }{\ }}%
\end{array}
& q^{8} & q^{7} & {q}^{5} & q^{5} & q^{7} & q^{6} & . & q^{5} & q^{4} & q^{3}%
& . & .\\
{\ }%
\begin{array}
[c]{c}%
{213\bar{2}{\ }}\\
{\bar{3}\bar{2}\bar{2}{\ \ }}\\
{\bar{2}\bar{1}{\ }{\ }{\ }}%
\end{array}
& q^{7} & q^{6} & q^{4} & q^{4} & q^{6} & q^{5} & . & q^{4} & q^{3} & q^{2} &
. & .\\
{\ }%
\begin{array}
[c]{c}%
{113\bar{2}{\ }}\\
{\bar{3}\bar{2}\bar{2}{\ \ }}\\
{\bar{1}\bar{1}{\ }{\ }{\ }}%
\end{array}
& q^{8} & q^{7} & q^{5} & q^{5} & {q}^{7} & q^{6} & . & q^{5} & q^{4} & q^{3}%
& . & .\\
{\ }%
\begin{array}
[c]{c}%
{313\bar{2}{\ }}\\
{\bar{3}\bar{3}\bar{2}{\ \ }}\\
{\bar{2}\bar{1}{\ }{\ }{\ }}%
\end{array}
& q^{8} & q^{7} & q^{5} & q^{5} & q^{7} & q^{6} & . & q^{5} & q^{4} & q^{3} &
. & .\\
{\ }%
\begin{array}
[c]{c}%
{113\bar{2}{\ }}\\
{\bar{2}\bar{3}\bar{2}{\ \ }}\\
{\bar{1}\bar{1}{\ }{\ }{\ }}%
\end{array}
& q^{9} & q^{8} & q^{6} & q^{6} & q^{8} & q^{7} & . & q^{6} & q^{5} & q^{4} &
. & .\\
{\ }%
\begin{array}
[c]{c}%
{313\bar{2}{\ }}\\
{\bar{3}\bar{3}\bar{2}{\ \ }}\\
{\bar{1}\bar{2}{\ }{\ }{\ }}%
\end{array}
& q^{9} & q^{8} & q^{6} & q^{6} & q^{8} & q^{7} & q^{4} & q^{6} & q^{5} &
q^{4} & . & .\\
{\ }%
\begin{array}
[c]{c}%
{213\bar{2}{\ }}\\
{\bar{2}\bar{3}\bar{2}{\ \ }}\\
{\bar{1}\bar{2}{\ }{\ }{\ }}%
\end{array}
& q^{10} & q^{9} & q^{7} & q^{7} & q^{9} & q^{8} & q^{5} & q^{7} & q^{6} &
q^{5} & . & .\\
{\ }%
\begin{array}
[c]{c}%
{\bar{3}13\bar{2}{\ }}\\
{\bar{2}3\bar{2}{\ \ }}\\
{\bar{1}\bar{3}{\ }{\ }{\ }}%
\end{array}
& . & . & . & . & . & . & q^{6} & . & . & . & . & .\\
{\ }%
\begin{array}
[c]{c}%
{\bar{3}13\bar{2}{\ }}\\
{\bar{2}2\bar{2}{\ \ }}\\
{\bar{1}\bar{2}{\ }{\ }{\ }}%
\end{array}
& . & . & . & . & . & . & q^{7} & . & . & . & . & .\\
{\ }%
\begin{array}
[c]{c}%
{1\bar{3}3\bar{2}{\ }}\\
{3\bar{2}\bar{3}{\ \ }}\\
{\bar{2}\bar{1}{\ }{\ }{\ }}%
\end{array}
& . & . & . & . & . & . & q^{4} & q^{2} & . & . & . & .
\end{array}
%
%
%
%
%
%
%
%
%
%
%
%
%
%
%
%
%
%
%
%
%
%
%
%
%
%
%
%
%
%
%
%
%
%
%
%
%
%
%
%
%
%
%
%
%
%
%
%
%
%
%
%
%
%
%
%
%
%
%
%
%
%
%
$}

\bigskip

{\scriptsize \hskip-16mm $
%
%
%
%
%
%
%
%
%
%
%
%
%
%
%
%
%
%
%
%
%
%
%
%
%
%
%
%
%
%
%
%
%
%
%
%
%
%
%
%
%
%
%
%
%
%
%
%
%
%
%
%
%
%
%
%
%
%
%
%
%
%
%
%
%
%
%
%
%
\begin{array}
[c]{ccccccccccccc}%
. & {\ }%
\begin{array}
[c]{c}%
{11\bar{2}\bar{1}}\text{ }\\
{3\bar{3}\bar{1}{\ \ }}\\
{\bar{2}\bar{2}{\ }{\ }{\ }}%
\end{array}
&
\begin{array}
[c]{c}%
{13\bar{3}\bar{1}}\text{ }\\
{3\bar{3}\bar{2}{\ \ }}\\
{\bar{2}\bar{2}{\ }{\ }{\ }}%
\end{array}
{\ } &
\begin{array}
[c]{c}%
{133\bar{1}}\text{ }\\
{\bar{3}\bar{3}\bar{2}{\ \ }}\\
{\bar{2}\bar{2}{\ }{\ }{\ }}%
\end{array}
{\ } &
\begin{array}
[c]{c}%
{13\bar{2}\bar{2}}\text{ }\\
{3\bar{3}\bar{1}{\ \ }}\\
{\bar{3}\bar{2}{\ }{\ }{\ }}%
\end{array}
{\ } & {\ }%
\begin{array}
[c]{c}%
{12\bar{2}\bar{2}}\text{ }\\
{3\bar{3}\bar{1}{\ \ }}\\
{\bar{2}\bar{2}{\ }{\ }{\ }}%
\end{array}
& {\ }%
\begin{array}
[c]{c}%
{13\bar{3}\bar{2}}\text{ }\\
{3\bar{3}\bar{1}{\ \ }}\\
{\bar{2}\bar{2}{\ }{\ }{\ }}%
\end{array}
&
\begin{array}
[c]{c}%
{13\bar{3}\bar{2}}\text{ }\\
{3\bar{2}\bar{2}{\ \ }}\\
{\bar{3}\bar{1}{\ }{\ }{\ }}%
\end{array}
{\ } & {\ }%
\begin{array}
[c]{c}%
{13\bar{3}\bar{2}}\text{ }\\
{3\bar{3}\bar{2}{\ \ }}\\
{\bar{2}\bar{1}{\ }{\ }{\ }}%
\end{array}
& {\ }%
\begin{array}
[c]{c}%
{133\bar{2}}\text{ }\\
{\bar{3}\bar{3}\bar{1}{\ \ }}\\
{\bar{2}\bar{2}{\ }{\ }{\ }}%
\end{array}
&
\begin{array}
[c]{c}%
{133\bar{2}}\text{ }\\
{\bar{3}\bar{3}\bar{2}{\ \ }}\\
{\bar{2}\bar{1}{\ }{\ }{\ }}%
\end{array}
{\ } & {\ }%
\begin{array}
[c]{c}%
{13\bar{3}\bar{3}}\text{ }\\
{3\bar{2}\bar{2}{\ \ }}\\
{\bar{2}\bar{1}{\ }{\ \ }}%
\end{array}
&
\begin{array}
[c]{c}%
{133\bar{3}}\text{ }\\
{\bar{3}\bar{2}\bar{2}{\ \ }}\\
{\bar{2}\bar{1}{\ }{\ }{\ }}%
\end{array}
{\ }\\
{\ }%
\begin{array}
[c]{c}%
{133\bar{2}{\ }}\\
{\bar{3}\bar{2}\bar{3}{\ \ }}\\
{\bar{2}\bar{1}{\ }{\ }{\ }}%
\end{array}
& . & . & q^{3} & q^{4} & . & . & q^{6} & q^{4} & . & q^{2} & . & .\\
{\ }%
\begin{array}
[c]{c}%
{333\bar{2}{\ }}\\
{\bar{3}\bar{3}\bar{3}{\ \ }}\\
{\bar{2}\bar{2}{\ }{\ }{\ }}%
\end{array}
& . & . & q^{4} & q^{5} & . & . & . & . & . & q^{3} & . & .\\
{\ }%
\begin{array}
[c]{c}%
{133\bar{2}{\ }}\\
{\bar{2}\bar{3}\bar{3}{\ \ }}\\
{\bar{1}\bar{2}{\ }{\ }{\ }}%
\end{array}
& . & . & q^{5} & q^{6} & . & . & . & . & . & q^{4} & . & .\\
{\ }%
\begin{array}
[c]{c}%
{313\bar{2}{\ }}\\
{\bar{3}\bar{2}\bar{3}{\ \ }}\\
{\bar{2}\bar{1}{\ }{\ }{\ }}%
\end{array}
& . & . & q^{5} & q^{6} & . & . & . & . & . & q^{4} & . & .\\
{\ }%
\begin{array}
[c]{c}%
{113\bar{2}{\ }}\\
{\bar{2}\bar{2}\bar{3}{\ \ }}\\
{\bar{1}\bar{1}{\ }{\ }{\ }}%
\end{array}
& {.} & . & q^{6} & q^{7} & . & . & . & . & . & q^{5} & . & .\\
{\ }%
\begin{array}
[c]{c}%
{313\bar{2}{\ }}\\
{\bar{2}\bar{3}\bar{3}{\ \ }}\\
{\bar{1}\bar{2}{\ }{\ }{\ }}%
\end{array}
& . & {.} & q^{7} & q^{8} & . & . & q^{6} & q^{4} & . & q^{6} & . & .\\
{\ }%
\begin{array}
[c]{c}%
{\bar{3}13\bar{2}{\ }}\\
{\bar{2}3\bar{3}{\ \ }}\\
{\bar{1}\bar{2}{\ }{\ }{\ }}%
\end{array}
& . & . & . & . & . & . & q^{8} & q^{6} & . & . & . & .\\
{\ }%
\begin{array}
[c]{c}%
{1\bar{3}2\bar{2}{\ }}\\
{3\bar{2}\bar{2}{\ \ }}\\
{\bar{2}\bar{1}{\ }{\ }{\ }}%
\end{array}
& . & q^{3} & {.} & q & . & q^{2} & q^{5} & q^{3} & . & . & . & .\\
{\ }%
\begin{array}
[c]{c}%
{132\bar{2}{\ }}\\
{\bar{3}\bar{2}\bar{2}{\ \ }}\\
{\bar{2}\bar{1}{\ }{\ }{\ }}%
\end{array}
& {.} & q^{5} & q^{3} & q^{5} & . & q^{4} & q^{7} & q^{5} & q^{2} & q^{3} &
. & .\\
{\ }%
\begin{array}
[c]{c}%
{332\bar{2}{\ }}\\
{\bar{3}\bar{3}\bar{2}{\ \ }}\\
{\bar{2}\bar{2}{\ }{\ }{\ }}%
\end{array}
& . & . & {q}^{4} & q^{6} & . & . & . & . & q^{3} & q^{4} & . & .\\
{\ }%
\begin{array}
[c]{c}%
{132\bar{2}{\ }}\\
{\bar{2}\bar{3}\bar{2}{\ \ }}\\
{\bar{1}\bar{2}{\ }{\ }{\ }}%
\end{array}
& {.} & . & q^{5} & q^{7} & . & . & . & . & q^{4} & q^{5} & . & .\\
{\ }%
\begin{array}
[c]{c}%
{312\bar{2}{\ }}\\
{\bar{3}\bar{2}\bar{2}{\ \ }}\\
{\bar{2}\bar{1}{\ }{\ }{\ }}%
\end{array}
& {.} & . & q^{5} & q^{7} & {.} & . & . & . & q^{4} & q^{5} & . & .\\
{\ }%
\begin{array}
[c]{c}%
{112\bar{2}{\ }}\\
{\bar{2}\bar{2}\bar{2}{\ \ }}\\
{\bar{1}\bar{1}{\ }{\ }{\ }}%
\end{array}
& . & . & q^{6} & q^{8} & . & . & . & . & q^{5} & q^{6} & . & .\\
{\ }%
\begin{array}
[c]{c}%
{312\bar{2}{\ }}\\
{\bar{2}\bar{3}\bar{2}{\ \ }}\\
{\bar{1}\bar{2}{\ }{\ }{\ }}%
\end{array}
& . & q^{5} & q^{7} & q^{9} & . & q^{4} & q^{7} & q^{5} & q^{6} & q^{7} & . &
.\\
{\ }%
\begin{array}
[c]{c}%
{\bar{3}12\bar{2}{\ }}\\
{\bar{2}3\bar{2}{\ \ }}\\
{\bar{1}\bar{2}{\ }{\ }{\ }}%
\end{array}
& . & q^{7} & . & {.} & . & q^{6} & q^{9} & q^{7} & . & . & . & .\\
{\ }%
\begin{array}
[c]{c}%
{1\bar{3}1\bar{2}{\ }}\\
{3\bar{2}\bar{1}{\ \ }}\\
{\bar{2}\bar{1}{\ }{\ }{\ }}%
\end{array}
& . & q^{4} & {.} & . & . & q^{3} & . & . & . & . & . & .\\
{\ }%
\begin{array}
[c]{c}%
{131\bar{2}{\ }}\\
{\bar{3}\bar{2}\bar{1}{\ \ }}\\
{\bar{2}\bar{1}{\ }{\ }{\ }}%
\end{array}
& {.} & q^{6} & q^{4} & . & . & q^{5} & . & . & . & . & . & .\\
{\ }%
\begin{array}
[c]{c}%
{331\bar{2}{\ }}\\
{\bar{3}\bar{3}\bar{1}{\ \ }}\\
{\bar{2}\bar{2}{\ }{\ }{\ }}%
\end{array}
& . & . & q^{5} & . & {.} & . & . & . & . & . & . & .\\
{\ }%
\begin{array}
[c]{c}%
{131\bar{2}{\ }}\\
{\bar{2}\bar{3}\bar{1}{\ \ }}\\
{\bar{1}\bar{2}{\ }{\ }{\ }}%
\end{array}
& {.} & . & {q}^{6} & . & . & . & . & . & . & . & . & .\\
{\ }%
\begin{array}
[c]{c}%
{311\bar{2}{\ }}\\
{\bar{3}\bar{2}\bar{1}{\ \ }}\\
{\bar{2}\bar{1}{\ }{\ }{\ }}%
\end{array}
& {.} & {.} & {q}^{6} & . & {.} & . & . & . & . & . & . & .\\
{\ }%
\begin{array}
[c]{c}%
{111\bar{2}{\ }}\\
{\bar{2}\bar{2}\bar{1}{\ \ }}\\
{\bar{1}\bar{1}{\ }{\ }{\ }}%
\end{array}
& . & . & q^{7} & . & . & . & . & . & . & . & . & .\\
{\ }%
\begin{array}
[c]{c}%
{311\bar{2}{\ }}\\
{\bar{2}\bar{3}\bar{1}{\ \ }}\\
{\bar{1}\bar{2}{\ }{\ }{\ }}%
\end{array}
& . & q^{6} & q^{8} & . & . & q^{5} & . & . & . & . & . & .\\
{\ }%
\begin{array}
[c]{c}%
{\bar{3}11\bar{2}{\ }}\\
{\bar{2}3\bar{1}{\ \ }}\\
{\bar{1}\bar{2}{\ }{\ }{\ }}%
\end{array}
& . & q^{8} & . & . & . & q^{7} & . & . & . & . & . & .\\
{\ }%
\begin{array}
[c]{c}%
{2\bar{3}1\bar{2}{\ }}\\
{3\bar{2}\bar{2}{\ \ }}\\
{\bar{2}\bar{1}{\ }{\ }{\ }}%
\end{array}
& . & q^{5} & . & . & . & q^{4} & . & q^{3} & . & . & . & .\\
{\ }%
\begin{array}
[c]{c}%
{1\bar{3}1\bar{2}{\ }}\\
{3\bar{2}\bar{2}{\ \ }}\\
{\bar{1}\bar{1}{\ }{\ }{\ }}%
\end{array}
& . & q^{6} & . & . & . & q^{5} & . & q^{4} & . & . & . & .
\end{array}
%
%
%
%
%
%
%
%
%
%
%
%
%
%
%
%
%
%
%
%
%
%
%
%
%
%
%
%
%
%
%
%
%
%
%
%
%
%
%
%
%
%
%
%
%
%
%
%
%
%
%
%
%
%
%
%
%
%
%
%
%
%
%
$}

\bigskip

{\scriptsize \hskip-16mm $
%
%
%
%
%
%
%
%
%
%
%
%
%
%
%
%
%
%
%
%
%
%
%
%
%
%
%
%
%
%
%
%
%
%
%
%
%
%
%
%
%
%
%
%
%
%
%
%
%
%
%
%
%
%
%
%
%
%
%
%
%
%
%
%
%
%
%
%
\begin{array}
[c]{ccccccccccccc}%
. & {\ }%
\begin{array}
[c]{c}%
{11\bar{2}\bar{1}}\text{ }\\
{3\bar{3}\bar{1}{\ \ }}\\
{\bar{2}\bar{2}{\ }{\ }{\ }}%
\end{array}
&
\begin{array}
[c]{c}%
{13\bar{3}\bar{1}}\text{ }\\
{3\bar{3}\bar{2}{\ \ }}\\
{\bar{2}\bar{2}{\ }{\ }{\ }}%
\end{array}
{\ } &
\begin{array}
[c]{c}%
{133\bar{1}}\text{ }\\
{\bar{3}\bar{3}\bar{2}{\ \ }}\\
{\bar{2}\bar{2}{\ }{\ }{\ }}%
\end{array}
{\ } &
\begin{array}
[c]{c}%
{13\bar{2}\bar{2}}\text{ }\\
{3\bar{3}\bar{1}{\ \ }}\\
{\bar{3}\bar{2}{\ }{\ }{\ }}%
\end{array}
{\ } & {\ }%
\begin{array}
[c]{c}%
{12\bar{2}\bar{2}}\text{ }\\
{3\bar{3}\bar{1}{\ \ }}\\
{\bar{2}\bar{2}{\ }{\ }{\ }}%
\end{array}
& {\ }%
\begin{array}
[c]{c}%
{13\bar{3}\bar{2}}\text{ }\\
{3\bar{3}\bar{1}{\ \ }}\\
{\bar{2}\bar{2}{\ }{\ }{\ }}%
\end{array}
&
\begin{array}
[c]{c}%
{13\bar{3}\bar{2}}\text{ }\\
{3\bar{2}\bar{2}{\ \ }}\\
{\bar{3}\bar{1}{\ }{\ }{\ }}%
\end{array}
{\ } & {\ }%
\begin{array}
[c]{c}%
{13\bar{3}\bar{2}}\text{ }\\
{3\bar{3}\bar{2}{\ \ }}\\
{\bar{2}\bar{1}{\ }{\ }{\ }}%
\end{array}
& {\ }%
\begin{array}
[c]{c}%
{133\bar{2}}\text{ }\\
{\bar{3}\bar{3}\bar{1}{\ \ }}\\
{\bar{2}\bar{2}{\ }{\ }{\ }}%
\end{array}
&
\begin{array}
[c]{c}%
{133\bar{2}}\text{ }\\
{\bar{3}\bar{3}\bar{2}{\ \ }}\\
{\bar{2}\bar{1}{\ }{\ }{\ }}%
\end{array}
{\ } & {\ }%
\begin{array}
[c]{c}%
{13\bar{3}\bar{3}}\text{ }\\
{3\bar{2}\bar{2}{\ \ }}\\
{\bar{2}\bar{1}{\ }{\ \ }}%
\end{array}
&
\begin{array}
[c]{c}%
{133\bar{3}}\text{ }\\
{\bar{3}\bar{2}\bar{2}{\ \ }}\\
{\bar{2}\bar{1}{\ }{\ }{\ }}%
\end{array}
{\ }\\
{\ }%
\begin{array}
[c]{c}%
{231\bar{2}{\ }}\\
{\bar{3}\bar{2}\bar{2}{\ \ }}\\
{\bar{2}\bar{1}{\ }{\ }{\ }}%
\end{array}
& . & q^{7} & q^{5} & . & . & q^{6} & . & q^{5} & q^{4} & q^{3} & . & .\\
{\ }%
\begin{array}
[c]{c}%
{131\bar{2}{\ }}\\
{\bar{3}\bar{2}\bar{2}{\ \ }}\\
{\bar{1}\bar{1}{\ }{\ }{\ }}%
\end{array}
& . & q^{8} & q^{6} & . & . & q^{7} & . & q^{6} & q^{5} & q^{4} & . & .\\
{\ }%
\begin{array}
[c]{c}%
{331\bar{2}{\ }}\\
{\bar{3}\bar{3}\bar{2}{\ \ }}\\
{\bar{2}\bar{1}{\ }{\ }{\ }}%
\end{array}
& . & . & q^{6} & . & . & . & . & . & q^{5} & q^{4} & . & .\\
{\ }%
\begin{array}
[c]{c}%
{131\bar{2}{\ }}\\
{\bar{2}\bar{3}\bar{2}{\ \ }}\\
{\bar{1}\bar{1}{\ }{\ }{\ }}%
\end{array}
& . & . & q^{7} & . & . & . & . & . & q^{6} & q^{5} & . & .\\
{\ }%
\begin{array}
[c]{c}%
{331\bar{2}{\ }}\\
{\bar{3}\bar{3}\bar{2}{\ \ }}\\
{\bar{1}\bar{2}{\ }{\ }{\ }}%
\end{array}
& {.} & . & q^{7} & . & . & . & . & . & q^{6} & q^{5} & . & .\\
{\ }%
\begin{array}
[c]{c}%
{231\bar{2}{\ }}\\
{\bar{2}\bar{3}\bar{2}{\ \ }}\\
{\bar{1}\bar{2}{\ }{\ }{\ }}%
\end{array}
& . & {.} & q^{8} & . & . & . & . & . & q^{7} & q^{6} & . & .\\
{\ }%
\begin{array}
[c]{c}%
{321\bar{2}{\ }}\\
{\bar{3}\bar{2}\bar{2}{\ \ }}\\
{\bar{2}\bar{1}{\ }{\ }{\ }}%
\end{array}
& . & . & q^{7} & . & . & . & . & . & q^{6} & q^{5} & . & .\\
{\ }%
\begin{array}
[c]{c}%
{121\bar{2}{\ }}\\
{\bar{2}\bar{2}\bar{2}{\ \ }}\\
{\bar{1}\bar{1}{\ }{\ }{\ }}%
\end{array}
& . & {.} & q^{8} & . & . & . & . & . & q^{7} & q^{6} & . & .\\
{\ }%
\begin{array}
[c]{c}%
{321\bar{2}{\ }}\\
{\bar{2}\bar{3}\bar{2}{\ \ }}\\
{\bar{1}\bar{2}{\ }{\ }{\ }}%
\end{array}
& {.} & q^{7} & q^{9} & . & . & q^{6} & . & q^{5} & q^{8} & q^{7} & . & .\\
{\ }%
\begin{array}
[c]{c}%
{\bar{3}21\bar{2}{\ }}\\
{\bar{2}3\bar{2}{\ \ }}\\
{\bar{1}\bar{2}{\ }{\ }{\ }}%
\end{array}
& . & q^{9} & . & . & . & q^{8} & . & q^{7} & . & . & . & .\\
{\ }%
\begin{array}
[c]{c}%
{311\bar{2}{\ }}\\
{\bar{3}\bar{2}\bar{2}{\ \ }}\\
{\bar{1}\bar{1}{\ }{\ }{\ }}%
\end{array}
& {.} & . & q^{8} & . & . & . & . & . & q^{7} & q^{6} & . & .\\
{\ }%
\begin{array}
[c]{c}%
{211\bar{2}{\ }}\\
{\bar{2}\bar{2}\bar{2}{\ \ }}\\
{\bar{1}\bar{1}{\ }{\ }{\ }}%
\end{array}
& {.} & . & q^{9} & . & {.} & . & . & . & q^{8} & q^{7} & . & .\\
{\ }%
\begin{array}
[c]{c}%
{311\bar{2}{\ }}\\
{\bar{2}\bar{3}\bar{2}{\ \ }}\\
{\bar{1}\bar{1}{\ }{\ }{\ }}%
\end{array}
& . & q^{8} & q^{10} & . & . & q^{7} & . & q^{6} & q^{9} & q^{8} & . & .\\
{\ }%
\begin{array}
[c]{c}%
{\bar{3}11\bar{2}{\ }}\\
{\bar{2}3\bar{2}{\ \ }}\\
{\bar{1}\bar{1}{\ }{\ }{\ }}%
\end{array}
& . & q^{10} & {.} & . & . & q^{9} & . & q^{8} & . & . & . & .\\
{\ }%
\begin{array}
[c]{c}%
{331\bar{2}{\ }}\\
{\bar{3}\bar{2}\bar{3}{\ \ }}\\
{\bar{2}\bar{1}{\ }{\ }{\ }}%
\end{array}
& . & . & . & {.} & . & . & . & q^{4} & . & q^{6} & . & .\\
{\ }%
\begin{array}
[c]{c}%
{131\bar{2}{\ }}\\
{\bar{2}\bar{2}\bar{3}{\ \ }}\\
{\bar{1}\bar{1}{\ }{\ }{\ }}%
\end{array}
& . & . & {.} & . & . & . & . & q^{5} & . & q^{7} & . & .\\
{\ }%
\begin{array}
[c]{c}%
{331\bar{2}{\ }}\\
{\bar{2}\bar{3}\bar{3}{\ \ }}\\
{\bar{1}\bar{2}{\ }{\ }{\ }}%
\end{array}
& {.} & . & {.} & {.} & . & . & . & q^{6} & . & q^{8} & . & .\\
{\ }%
\begin{array}
[c]{c}%
{311\bar{2}{\ }}\\
{\bar{2}\bar{2}\bar{3}{\ \ }}\\
{\bar{1}\bar{1}{\ }{\ }{\ }}%
\end{array}
& . & . & . & {.} & {.} & . & . & q^{7} & . & q^{9} & . & .\\
{\ }%
\begin{array}
[c]{c}%
{3\bar{3}1\bar{2}{\ }}\\
{\bar{3}\bar{2}3{\ \ }}\\
{\bar{2}\bar{1}{\ }{\ }{\ }}%
\end{array}
& {.} & . & {.} & {.} & . & . & . & q^{6} & . & . & . & .\\
{\ }%
\begin{array}
[c]{c}%
{1\bar{3}1\bar{2}{\ }}\\
{\bar{2}\bar{2}3{\ \ }}\\
{\bar{1}\bar{1}{\ }{\ }{\ }}%
\end{array}
& {.} & {.} & {.} & {.} & {.} & . & . & q^{7} & . & . & . & .\\
{\ }%
\begin{array}
[c]{c}%
{\bar{3}31\bar{2}{\ }}\\
{\bar{2}\bar{3}3{\ \ }}\\
{\bar{1}\bar{2}{\ }{\ }{\ }}%
\end{array}
& . & . & . & . & . & . & . & q^{8} & . & . & . & .\\
{\ }%
\begin{array}
[c]{c}%
{\bar{3}11\bar{2}{\ }}\\
{\bar{2}\bar{2}3{\ \ }}\\
{\bar{1}\bar{1}{\ }{\ }{\ }}%
\end{array}
& . & . & . & . & . & . & . & q^{9} & . & . & . & .\\
{\ }%
\begin{array}
[c]{c}%
{13\bar{2}\bar{3}{\ }}\\
{3\bar{3}\bar{1}{\ \ }}\\
{\bar{2}\bar{2}{\ }{\ }{\ }}%
\end{array}
& . & . & . & q^{4} & q^{2} & q & . & . & . & . & . & .\\
{\ }%
\begin{array}
[c]{c}%
{11\bar{2}\bar{3}{\ }}\\
{3\bar{2}1{\ \ }}\\
{\bar{2}\bar{1}{\ }{\ }{\ }}%
\end{array}
& . & . & . & q^{5} & q^{3} & q^{2} & . & . & . & . & . & .\\
{\ }%
\begin{array}
[c]{c}%
{31\bar{2}\bar{3}{\ }}\\
{\bar{3}3\bar{1}{\ \ }}\\
{\bar{2}\bar{2}{\ }{\ }{\ }}%
\end{array}
& . & . & . & q^{6} & q^{4} & q^{3} & . & . & . & . & . & .
\end{array}
%
%
%
%
%
%
%
%
%
%
%
%
%
%
%
%
%
%
%
%
%
%
%
%
%
%
%
%
%
%
%
%
%
%
%
%
%
%
%
%
%
%
%
%
%
%
%
%
%
%
%
%
%
%
%
%
%
%
%
%
%
%
%
$}

\bigskip

{\scriptsize \hskip-16mm $
%
%
%
%
%
%
%
%
%
%
%
%
%
%
%
%
%
%
%
%
%
%
%
%
%
%
%
%
%
%
%
%
%
%
%
%
%
%
%
%
%
%
%
%
%
%
%
%
%
%
%
%
%
%
%
%
%
%
%
%
%
%
%
%
%
%
%
%
\begin{array}
[c]{ccccccccccccc}%
. & {\ }%
\begin{array}
[c]{c}%
{11\bar{2}\bar{1}}\text{ }\\
{3\bar{3}\bar{1}{\ \ }}\\
{\bar{2}\bar{2}{\ }{\ }{\ }}%
\end{array}
&
\begin{array}
[c]{c}%
{13\bar{3}\bar{1}}\text{ }\\
{3\bar{3}\bar{2}{\ \ }}\\
{\bar{2}\bar{2}{\ }{\ }{\ }}%
\end{array}
{\ } &
\begin{array}
[c]{c}%
{133\bar{1}}\text{ }\\
{\bar{3}\bar{3}\bar{2}{\ \ }}\\
{\bar{2}\bar{2}{\ }{\ }{\ }}%
\end{array}
{\ } &
\begin{array}
[c]{c}%
{13\bar{2}\bar{2}}\text{ }\\
{3\bar{3}\bar{1}{\ \ }}\\
{\bar{3}\bar{2}{\ }{\ }{\ }}%
\end{array}
{\ } & {\ }%
\begin{array}
[c]{c}%
{12\bar{2}\bar{2}}\text{ }\\
{3\bar{3}\bar{1}{\ \ }}\\
{\bar{2}\bar{2}{\ }{\ }{\ }}%
\end{array}
& {\ }%
\begin{array}
[c]{c}%
{13\bar{3}\bar{2}}\text{ }\\
{3\bar{3}\bar{1}{\ \ }}\\
{\bar{2}\bar{2}{\ }{\ }{\ }}%
\end{array}
&
\begin{array}
[c]{c}%
{13\bar{3}\bar{2}}\text{ }\\
{3\bar{2}\bar{2}{\ \ }}\\
{\bar{3}\bar{1}{\ }{\ }{\ }}%
\end{array}
{\ } & {\ }%
\begin{array}
[c]{c}%
{13\bar{3}\bar{2}}\text{ }\\
{3\bar{3}\bar{2}{\ \ }}\\
{\bar{2}\bar{1}{\ }{\ }{\ }}%
\end{array}
& {\ }%
\begin{array}
[c]{c}%
{133\bar{2}}\text{ }\\
{\bar{3}\bar{3}\bar{1}{\ \ }}\\
{\bar{2}\bar{2}{\ }{\ }{\ }}%
\end{array}
&
\begin{array}
[c]{c}%
{133\bar{2}}\text{ }\\
{\bar{3}\bar{3}\bar{2}{\ \ }}\\
{\bar{2}\bar{1}{\ }{\ }{\ }}%
\end{array}
{\ } & {\ }%
\begin{array}
[c]{c}%
{13\bar{3}\bar{3}}\text{ }\\
{3\bar{2}\bar{2}{\ \ }}\\
{\bar{2}\bar{1}{\ }{\ \ }}%
\end{array}
&
\begin{array}
[c]{c}%
{133\bar{3}}\text{ }\\
{\bar{3}\bar{2}\bar{2}{\ \ }}\\
{\bar{2}\bar{1}{\ }{\ }{\ }}%
\end{array}
{\ }\\
{\ }%
\begin{array}
[c]{c}%
{11\bar{2}\bar{3}{\ }}\\
{\bar{2}3\bar{1}{\ \ }}\\
{\bar{1}\bar{2}{\ }{\ }{\ }}%
\end{array}
& . & . & . & q^{7} & q^{5} & q^{4} & . & . & . & . & . & .\\
{\ }%
\begin{array}
[c]{c}%
\mathbf{13\bar{3}\bar{3}}{{\ }}\\
\mathbf{3\bar{2}\bar{2}\ \ }\\
\mathbf{\bar{2}\bar{1}\ \ \ }%
\end{array}
& . & . & . & q^{6} & q^{4} & q+q^{3} & q^{4} & q^{2} & . & . & 1 & .\\
{\ }%
\begin{array}
[c]{c}%
{31\bar{3}\bar{3}{\ }}\\
{\bar{2}3\bar{2}{\ \ }}\\
{\bar{1}\bar{2}{\ }{\ }{\ }}%
\end{array}
& . & . & . & q^{8} & q^{6} & q^{3}+q^{5} & q^{6} & q^{4} & . & . & q^{2} &
.\\
{\ }%
\begin{array}
[c]{c}%
{1\bar{3}3\bar{3}{\ }}\\
{3\bar{2}\bar{2}{\ \ }}\\
{\bar{2}\bar{1}{\ }{\ }{\ }}%
\end{array}
& . & . & . & . & . & q^{3} & q^{6} & q^{4} & . & . & q^{2} & .\\
{\ }%
\begin{array}
[c]{c}%
\mathbf{133\bar{3}}{{\ }}\\
\mathbf{\bar{3}\bar{2}\bar{2}\ \ }\\
\mathbf{\bar{2}\bar{1}\ \ \ }%
\end{array}
& {.} & . & . & q^{6} & q^{4} & q+q^{5}+q^{3} & q^{8} & q^{2}+q^{6} & q^{3} &
q^{4} & q^{4} & 1\\
{\ }%
\begin{array}
[c]{c}%
{333\bar{3}{\ }}\\
{\bar{3}\bar{3}\bar{2}{\ \ }}\\
{\bar{2}\bar{2}{\ }{\ }{\ }}%
\end{array}
& . & {.} & . & q^{7} & q^{5} & q^{2}+q^{4} & . & q^{3} & q^{4} & q^{5} & . &
q\\
{\ }%
\begin{array}
[c]{c}%
{133\bar{3}{\ }}\\
{\bar{2}\bar{3}\bar{2}{\ \ }}\\
{\bar{1}\bar{2}{\ }{\ }{\ }}%
\end{array}
& . & . & . & q^{8} & q^{6} & q^{3}+q^{5} & . & q^{4} & q^{5} & q^{6} & . &
q^{2}\\
{\ }%
\begin{array}
[c]{c}%
{313\bar{3}{\ }}\\
{\bar{3}\bar{2}\bar{2}{\ \ }}\\
{\bar{2}\bar{1}{\ }{\ }{\ }}%
\end{array}
& . & {.} & {.} & q^{8} & q^{6} & q^{3}+q^{5} & . & q^{4} & q^{5} & q^{6} &
. & q^{2}\\
{\ }%
\begin{array}
[c]{c}%
{113\bar{3}{\ }}\\
{\bar{2}\bar{2}\bar{2}{\ \ }}\\
{\bar{1}\bar{1}{\ }{\ }{\ }}%
\end{array}
& {.} & {.} & . & q^{9} & q^{7} & q^{4}+q^{6} & . & q^{5} & q^{6} & q^{7} &
. & q^{3}\\
{\ }%
\begin{array}
[c]{c}%
{313\bar{3}{\ }}\\
{\bar{2}\bar{3}\bar{2}{\ \ }}\\
{\bar{1}\bar{2}{\ }{\ }{\ }}%
\end{array}
& . & . & {.} & q^{10} & q^{8} & 2q^{5}+q^{7} & q^{8} & 2q^{6} & q^{7} & q^{8}%
& q^{4} & q^{4}\\
{\ }%
\begin{array}
[c]{c}%
{\bar{3}13\bar{3}{\ }}\\
{\bar{2}3\bar{2}{\ \ }}\\
{\bar{1}\bar{2}{\ }{\ }{\ }}%
\end{array}
& {.} & . & {.} & {.} & . & q^{7} & q^{10} & q^{8} & . & . & q^{6} & .\\
{\ }%
\begin{array}
[c]{c}%
{331\bar{3}{\ }}\\
{\bar{3}\bar{2}\bar{2}{\ \ }}\\
{\bar{2}\bar{1}{\ }{\ }{\ }}%
\end{array}
& {.} & . & . & . & {.} & q^{4} & . & q^{5} & q^{6} & q^{7} & . & q^{3}\\
{\ }%
\begin{array}
[c]{c}%
{131\bar{3}{\ }}\\
{\bar{2}\bar{2}\bar{2}{\ \ }}\\
{\bar{1}\bar{1}{\ }{\ }{\ }}%
\end{array}
& . & . & . & . & . & q^{5} & . & q^{6} & q^{7} & q^{8} & . & q^{4}\\
{\ }%
\begin{array}
[c]{c}%
{331\bar{3}{\ }}\\
{\bar{2}\bar{3}\bar{2}{\ \ }}\\
{\bar{1}\bar{2}{\ }{\ }{\ }}%
\end{array}
& . & . & {.} & . & . & q^{6} & . & q^{7} & q^{8} & q^{9} & . & q^{5}\\
{\ }%
\begin{array}
[c]{c}%
{311\bar{3}{\ }}\\
{\bar{2}\bar{2}\bar{2}{\ \ }}\\
{\bar{1}\bar{1}{\ }{\ }{\ }}%
\end{array}
& . & . & . & . & . & q^{7} & . & q^{8} & q^{9} & q^{10} & . & q^{6}\\
{\ }%
\begin{array}
[c]{c}%
{13\bar{2}3{\ }}\\
{\bar{3}\bar{3}\bar{1}{\ \ }}\\
{\bar{2}\bar{2}{\ }{\ }{\ }}%
\end{array}
& . & . & {.} & . & q^{4} & q^{3} & . & . & q & . & . & .\\
{\ }%
\begin{array}
[c]{c}%
{11\bar{2}3{\ }}\\
{\bar{3}\bar{2}\bar{1}{\ \ }}\\
{\bar{2}\bar{1}{\ }{\ }{\ }}%
\end{array}
& {.} & . & {.} & . & q^{5} & q^{4} & . & . & q^{2} & . & . & .\\
{\ }%
\begin{array}
[c]{c}%
{31\bar{2}3{\ }}\\
{\bar{3}\bar{3}\bar{1}{\ \ }}\\
{\bar{2}\bar{2}{\ }{\ }{\ }}%
\end{array}
& . & . & . & . & q^{6} & q^{5} & . & . & q^{3} & . & . & .\\
{\ }%
\begin{array}
[c]{c}%
{11\bar{2}3{\ }}\\
{\bar{2}\bar{3}\bar{1}{\ \ }}\\
{\bar{1}\bar{2}{\ }{\ }{\ }}%
\end{array}
& {.} & . & {.} & . & q^{7} & q^{6} & . & . & q^{4} & . & . & .\\
{\ }%
\begin{array}
[c]{c}%
{1\bar{3}\bar{3}3{\ }}\\
{3\bar{2}\bar{2}{\ \ }}\\
{\bar{2}\bar{1}{\ }{\ }{\ }}%
\end{array}
& {.} & {.} & {.} & . & . & q^{3} & . & . & . & . & q^{4} & .\\
{\ }%
\begin{array}
[c]{c}%
13\bar{3}3{{\ }}\\
\bar{3}\bar{2}\bar{2}\ \ \\
\bar{2}\bar{1}\ \ \
\end{array}
& . & . & . & . & q^{6} & q^{3}+2q^{5} & . & q^{4} & q^{3} & . & q^{6} &
q^{2}\\
{\ }%
\begin{array}
[c]{c}%
{33\bar{3}3{\ }}\\
{\bar{3}\bar{3}\bar{2}{\ \ }}\\
{\bar{2}\bar{2}{\ }{\ }{\ }}%
\end{array}
& . & . & . & . & q^{7} & q^{4}+q^{6} & . & q^{5} & q^{4} & . & . & q^{3}\\
{\ }%
\begin{array}
[c]{c}%
{13\bar{3}3{\ }}\\
{\bar{2}\bar{3}\bar{2}{\ \ }}\\
{\bar{1}\bar{2}{\ }{\ }{\ }}%
\end{array}
& . & . & . & . & q^{8} & q^{5}+q^{7} & . & q^{6} & q^{5} & . & . & q^{4}\\
{\ }%
\begin{array}
[c]{c}%
{31\bar{3}3{\ }}\\
{\bar{3}\bar{2}\bar{2}{\ \ }}\\
{\bar{2}\bar{1}{\ }{\ }{\ }}%
\end{array}
& . & . & . & . & q^{8} & q^{5}+q^{7} & . & q^{6} & q^{5} & . & . & q^{4}\\
{\ }%
\begin{array}
[c]{c}%
{11\bar{3}3{\ }}\\
{\bar{2}\bar{2}\bar{2}{\ \ }}\\
{\bar{1}\bar{1}{\ }{\ }{\ }}%
\end{array}
& . & . & . & . & q^{9} & q^{6}+q^{8} & . & q^{7} & q^{6} & . & . & q^{5}%
\end{array}
%
%
%
%
%
%
%
%
%
%
%
%
%
%
%
%
%
%
%
%
%
%
%
%
%
%
%
%
%
%
%
%
%
%
%
%
%
%
%
%
%
%
%
%
%
%
%
%
%
%
%
%
%
%
%
%
%
%
%
%
%
%
%
$}

\bigskip

{\scriptsize \hskip-16mm $
%
%
%
%
%
%
%
%
%
%
%
%
%
%
%
%
%
%
%
%
%
%
%
%
%
%
%
%
%
%
%
%
%
%
%
%
%
%
%
%
%
%
%
%
%
%
%
%
%
%
%
%
%
%
%
%
%
%
%
%
%
%
%
%
%
%
%
%
\begin{array}
[c]{ccccccccccccc}%
. & {\ }%
\begin{array}
[c]{c}%
{11\bar{2}\bar{1}}\text{ }\\
{3\bar{3}\bar{1}{\ \ }}\\
{\bar{2}\bar{2}{\ }{\ }{\ }}%
\end{array}
&
\begin{array}
[c]{c}%
{13\bar{3}\bar{1}}\text{ }\\
{3\bar{3}\bar{2}{\ \ }}\\
{\bar{2}\bar{2}{\ }{\ }{\ }}%
\end{array}
{\ } &
\begin{array}
[c]{c}%
{133\bar{1}}\text{ }\\
{\bar{3}\bar{3}\bar{2}{\ \ }}\\
{\bar{2}\bar{2}{\ }{\ }{\ }}%
\end{array}
{\ } &
\begin{array}
[c]{c}%
{13\bar{2}\bar{2}}\text{ }\\
{3\bar{3}\bar{1}{\ \ }}\\
{\bar{3}\bar{2}{\ }{\ }{\ }}%
\end{array}
{\ } & {\ }%
\begin{array}
[c]{c}%
{12\bar{2}\bar{2}}\text{ }\\
{3\bar{3}\bar{1}{\ \ }}\\
{\bar{2}\bar{2}{\ }{\ }{\ }}%
\end{array}
& {\ }%
\begin{array}
[c]{c}%
{13\bar{3}\bar{2}}\text{ }\\
{3\bar{3}\bar{1}{\ \ }}\\
{\bar{2}\bar{2}{\ }{\ }{\ }}%
\end{array}
&
\begin{array}
[c]{c}%
{13\bar{3}\bar{2}}\text{ }\\
{3\bar{2}\bar{2}{\ \ }}\\
{\bar{3}\bar{1}{\ }{\ }{\ }}%
\end{array}
{\ } & {\ }%
\begin{array}
[c]{c}%
{13\bar{3}\bar{2}}\text{ }\\
{3\bar{3}\bar{2}{\ \ }}\\
{\bar{2}\bar{1}{\ }{\ }{\ }}%
\end{array}
& {\ }%
\begin{array}
[c]{c}%
{133\bar{2}}\text{ }\\
{\bar{3}\bar{3}\bar{1}{\ \ }}\\
{\bar{2}\bar{2}{\ }{\ }{\ }}%
\end{array}
&
\begin{array}
[c]{c}%
{133\bar{2}}\text{ }\\
{\bar{3}\bar{3}\bar{2}{\ \ }}\\
{\bar{2}\bar{1}{\ }{\ }{\ }}%
\end{array}
{\ } & {\ }%
\begin{array}
[c]{c}%
{13\bar{3}\bar{3}}\text{ }\\
{3\bar{2}\bar{2}{\ \ }}\\
{\bar{2}\bar{1}{\ }{\ \ }}%
\end{array}
&
\begin{array}
[c]{c}%
{133\bar{3}}\text{ }\\
{\bar{3}\bar{2}\bar{2}{\ \ }}\\
{\bar{2}\bar{1}{\ }{\ }{\ }}%
\end{array}
{\ }\\
{\ }%
\begin{array}
[c]{c}%
{31\bar{3}3{\ }}\\
{\bar{2}\bar{3}\bar{2}{\ \ }}\\
{\bar{1}\bar{2}{\ }{\ }{\ }}%
\end{array}
& . & . & . & . & q^{10} & q^{5}+q^{7}+q^{9} & . & q^{8} & q^{7} & . & q^{6} &
q^{6}\\
{\ }%
\begin{array}
[c]{c}%
{\bar{3}1\bar{3}3{\ }}\\
{\bar{2}3\bar{2}{\ \ }}\\
{\bar{1}\bar{2}{\ }{\ }{\ }}%
\end{array}
& . & . & . & . & . & q^{7} & . & . & . & . & q^{8} & .\\
{\ }%
\begin{array}
[c]{c}%
{1\bar{3}33{\ }}\\
{\bar{3}\bar{2}\bar{2}{\ \ }}\\
{\bar{2}\bar{1}{\ }{\ }{\ }}%
\end{array}
& . & . & . & . & . & q^{5}+q^{7} & . & q^{6} & . & . & q^{8} & q^{4}\\
{\ }%
\begin{array}
[c]{c}%
{\bar{3}133{\ }}\\
{\bar{2}\bar{3}\bar{2}{\ \ }}\\
{\bar{1}\bar{2}{\ }{\ }{\ }}%
\end{array}
& . & . & . & . & . & q^{7}+q^{9} & . & q^{8} & . & . & q^{10} & q^{6}\\
{\ }%
\begin{array}
[c]{c}%
{3\bar{3}13{\ }}\\
{\bar{3}\bar{2}\bar{2}{\ \ }}\\
{\bar{2}1{\ }{\ }{\ }}%
\end{array}
& {.} & . & . & . & . & q^{6} & . & q^{7} & . & . & . & q^{5}\\
{\ }%
\begin{array}
[c]{c}%
1{\bar{3}13{\ }}\\
{\bar{2}\bar{2}\bar{2}{\ \ }}\\
{\bar{1}\bar{1}{\ }{\ }{\ }}%
\end{array}
& . & {.} & . & . & . & q^{7} & . & q^{8} & . & . & . & q^{6}\\
{\ }%
\begin{array}
[c]{c}%
{\bar{3}313{\ }}\\
{\bar{2}\bar{3}\bar{2}{\ \ }}\\
{\bar{1}\bar{2}{\ }{\ }{\ }}%
\end{array}
& . & . & . & . & . & q^{8} & . & q^{9} & . & . & . & q^{7}\\
{\ }%
\begin{array}
[c]{c}%
{\bar{3}113{\ }}\\
{\bar{2}\bar{2}\bar{2}{\ \ }}\\
{\bar{1}\bar{1}{\ }{\ }{\ }}%
\end{array}
& . & {.} & {.} & . & . & q^{9} & . & q^{10} & . & . & . & q^{8}\\
{\ }%
\begin{array}
[c]{c}%
{1\bar{3}\bar{2}1{\ }}\\
{3\bar{2}\bar{1}{\ \ }}\\
{\bar{2}\bar{1}{\ }{\ }{\ }}%
\end{array}
& {.} & {.} & . & {.} & . & q^{4} & . & . & . & . & . & .\\
{\ }%
\begin{array}
[c]{c}%
{13\bar{2}1{\ }}\\
{\bar{3}\bar{2}\bar{1}{\ \ }}\\
{\bar{2}\bar{1}{\ }{\ }{\ }}%
\end{array}
& . & . & {.} & . & . & q^{6} & . & . & q^{4} & . & . & .\\
{\ }%
\begin{array}
[c]{c}%
{33\bar{2}1{\ }}\\
{\bar{3}\bar{3}\bar{1}{\ \ }}\\
{\bar{2}\bar{2}{\ }{\ }{\ }}%
\end{array}
& {.} & . & {.} & . & . & . & . & . & q^{5} & . & . & .\\
{\ }%
\begin{array}
[c]{c}%
{13\bar{2}1{\ }}\\
{\bar{2}\bar{3}\bar{1}{\ \ }}\\
{\bar{1}\bar{2}{\ }{\ }{\ }}%
\end{array}
& {.} & . & . & . & {.} & . & . & . & q^{6} & . & . & .\\
{\ }%
\begin{array}
[c]{c}%
{31\bar{2}1{\ }}\\
{\bar{3}\bar{2}\bar{1}{\ \ }}\\
{\bar{2}\bar{1}{\ }{\ }{\ }}%
\end{array}
& . & . & . & . & . & . & . & . & q^{6} & . & . & .\\
{\ }%
\begin{array}
[c]{c}%
{11\bar{2}1{\ }}\\
{\bar{2}\bar{2}1{\ \ }}\\
{\bar{1}\bar{1}{\ }{\ }{\ }}%
\end{array}
& . & . & {.} & . & . & . & . & . & q^{7} & . & . & .\\
{\ }%
\begin{array}
[c]{c}%
{31\bar{2}1{\ }}\\
{\bar{2}\bar{3}\bar{1}{\ \ }}\\
{\bar{1}\bar{2}{\ }{\ }{\ }}%
\end{array}
& . & . & . & . & . & q^{6} & . & . & q^{8} & . & . & .\\
{\ }%
\begin{array}
[c]{c}%
{\bar{3}1\bar{2}1{\ }}\\
{\bar{2}3\bar{1}{\ \ }}\\
{\bar{1}\bar{2}{\ }{\ }{\ }}%
\end{array}
& . & . & {.} & . & . & q^{8} & . & . & . & . & . & .\\
{\ }%
\begin{array}
[c]{c}%
{33\bar{3}1{\ }}\\
{\bar{3}\bar{2}\bar{2}{\ \ }}\\
{\bar{2}\bar{1}{\ }{\ }{\ }}%
\end{array}
& {.} & . & {.} & . & . & q^{5} & . & . & q^{7} & . & . & q^{4}\\
{\ }%
\begin{array}
[c]{c}%
{13\bar{3}1{\ }}\\
{\bar{2}\bar{2}\bar{2}{\ \ }}\\
{\bar{1}\bar{1}{\ }{\ }{\ }}%
\end{array}
& . & . & . & {.} & {.} & q^{6} & . & . & q^{8} & . & . & q^{5}\\
{\ }%
\begin{array}
[c]{c}%
{33\bar{3}1{\ }}\\
{\bar{2}\bar{3}\bar{2}{\ \ }}\\
{\bar{1}\bar{2}{\ }{\ }{\ }}%
\end{array}
& {.} & . & {.} & {.} & . & q^{7} & . & . & q^{9} & . & . & q^{6}\\
{\ }%
\begin{array}
[c]{c}%
{31\bar{3}1{\ }}\\
{\bar{2}\bar{2}\bar{2}{\ \ }}\\
{\bar{1}\bar{1}{\ }{\ }{\ }}%
\end{array}
& {.} & {.} & {.} & {.} & {.} & q^{8} & . & . & q^{10} & . & . & q^{7}\\
{\ }%
\begin{array}
[c]{c}%
{3\bar{3}31{\ }}\\
{\bar{3}\bar{2}\bar{2}{\ \ }}\\
{\bar{2}\bar{1}{\ }{\ }{\ }}%
\end{array}
& . & . & . & . & . & q^{7} & . & . & . & . & . & q^{6}\\
{\ }%
\begin{array}
[c]{c}%
{1\bar{3}31{\ }}\\
{\bar{2}\bar{2}\bar{2}{\ \ }}\\
{\bar{1}\bar{1}{\ }{\ }{\ }}%
\end{array}
& . & . & . & . & . & q^{8} & . & . & . & . & . & q^{7}\\
{\ }%
\begin{array}
[c]{c}%
{\bar{3}331{\ }}\\
{\bar{2}\bar{3}\bar{2}{\ \ }}\\
{\bar{1}\bar{2}{\ }{\ }{\ }}%
\end{array}
& . & . & . & . & . & q^{9} & . & . & . & . & . & q^{8}\\
{\ }%
\begin{array}
[c]{c}%
{\bar{3}331{\ }}\\
{\bar{2}\bar{2}\bar{2}{\ \ }}\\
{\bar{1}\bar{1}{\ }{\ }{\ }}%
\end{array}
& . & . & . & . & . & q^{10} & . & . & . & . & . & q^{9}\\
{\ }%
\begin{array}
[c]{c}%
{3\bar{3}11{\ }}\\
{\bar{2}\bar{2}\bar{2}{\ \ }}\\
{\bar{1}\bar{1}{\ }{\ }{\ }}%
\end{array}
& . & . & . & . & . & . & . & . & . & . & . & q^{8}\\
{\ }%
\begin{array}
[c]{c}%
{\bar{3}311{\ }}\\
{\bar{2}\bar{2}\bar{2}{\ \ }}\\
{\bar{1}\bar{1}{\ }{\ }{\ }}%
\end{array}
& . & . & . & . & . & . & . & . & . & . & . & q^{10}%
\end{array}
%
%
%
%
%
%
%
%
%
%
%
%
%
%
%
%
%
%
%
%
%
%
%
%
%
%
%
%
%
%
%
%
%
%
%
%
%
%
%
%
%
%
%
%
%
%
%
%
%
%
%
%
%
%
%
%
%
%
%
%
%
%
%
$}

\bigskip

\noindent Note that all the coefficients $d_{\tau,T}(q)$ of the above matrix
are in $\mathbb{N}[q].\;$This not true in general. For example, consider the
canonical basis of the weight space of the $U_{q}(sp_{8})$-module $V(1,1,1,1)$
corresponding to the weight $(0,0,0,0)$.\ Then for $T=%
\begin{tabular}
[c]{|l|l|}\hline
$\mathtt{1}$ & $\mathtt{4}$\\\hline
$\mathtt{3}$ & $\mathtt{\bar{4}}$\\\hline
$\mathtt{4}$ & $\mathtt{\bar{3}}$\\\hline
$\mathtt{\bar{4}}$ & $\mathtt{\bar{1}}$\\\hline
\end{tabular}
$, there are two coefficients $d_{\tau,T}(q)\notin\mathbb{N[}q]$ in $G(T).$
More precisely, we have $d_{\tau_{1},T}(q)=-q^{4}$ and $d_{\tau_{2}%
,T}(q)=-q^{4}$ for $\tau_{1}=%
\begin{tabular}
[c]{|l|l|}\hline
$\mathtt{1}$ & $\mathtt{2}$\\\hline
$\mathtt{4}$ & $\mathtt{3}$\\\hline
$\mathtt{\bar{3}}$ & $\mathtt{\bar{4}}$\\\hline
$\mathtt{\bar{2}}$ & $\mathtt{\bar{1}}$\\\hline
\end{tabular}
$ and $\tau_{2}=%
\begin{tabular}
[c]{|l|l|}\hline
$\mathtt{2}$ & $\mathtt{1}$\\\hline
$\mathtt{3}$ & $\mathtt{4}$\\\hline
$\mathtt{\bar{4}}$ & $\mathtt{\bar{3}}$\\\hline
$\mathtt{\bar{1}}$ & $\mathtt{\bar{2}}$\\\hline
\end{tabular}
.$

\noindent\textbf{Acknowledgments:} We are very grateful to P. Toffin who
implemented our algorithm in AXIOM and pointed out the above example of
non-positive coefficients.`


\begin{thebibliography}{99}
\bibitem{Ch-Pr}\textsc{V. Chari, A. Presley}, \textit{A guide to quantum
groups}, Cambridge University Press 1994.

\bibitem {Jan}\textsc{J. C. Jantzen}, \textit{Lectures on quantum groups},
Graduate Studies in Math. \textbf{6}, A.M.S 1995

\bibitem {Ka0}\textsc{M. Kashiwara}, \textit{Crystallizing the }%
$q$\textit{-analogue of universal enveloping algebra}, Commun. Math. Phys,
\textbf{133} (1990), 249-260.

\bibitem {Ka}\textsc{M. Kashiwara}, \textit{On crystal bases of the }%
$q$\textit{-analogue of universal enveloping algebras}, Duke Math. J,
\textbf{63} (1991), 465-516.

\bibitem {Ka1}\textsc{M. Kashiwara,} \textit{Crystallization of quantized
universal enveloping algebras}, Sugaku Expositiones, \textbf{7} (1994), 99-115

\bibitem {Ka2}\textsc{M. Kashiwara,}\textit{ On crystal bases, }Canadian
Mathematical Society, Conference Proceedings, \textbf{16} (1995), 155-197.

\bibitem {KN}\textsc{M. Kashiwara, T. Nakashima,} \textit{Crystal graphs for
representations of the }$q$\textit{-analogue of classical Lie algebras},
Journal of Algebra, \textbf{165} (1994), 295-345.

\bibitem {L-T}\textsc{B. Leclerc, P. Toffin}, \textit{A simple algorithm for
computing the global crystal basis of an irreducible }$U_{q}(sl_{n}%
)$\textit{-module, }Int. J.\ Algebra Computation, \textbf{10} (2000), 191-208.

\bibitem {Lec}\textsc{C. Lecouvey,} \textit{Schensted-type correspondence,
Plactic Monoid and Jeu de Taquin for type }$C_{n},$ Preprint 1999.

\bibitem {Lut}\textsc{G. Lusztig}, \textit{Quivers, perverse sheaves, and
quantized enveloping algebras, }J. Am. Math. Soc, \textbf{4} (1991), 365-421.

\bibitem {Ma}\textsc{R. Marsh}, \textit{Algorithms to obtain the canonical
basis in some fundamental modules of quantum groups}, Journal of Algebra
\textbf{196,} 831-860 (1996)

\bibitem {Sh}\textsc{J.T. Sheats,} \textit{A symplectic Jeu de Taquin
bijection between the tableaux of King and De Concini}, Trans. A.M.S,
\textbf{351} (1999), 3569-3607.

\bibitem {ZR}\textsc{A.V. Zelevinsky, V.S. Retakh, }\textit{The base affine
space and canonical basis in irreducible representations of the group }%
$Sp_{4},$ Sov. Math, \textbf{37} (1988), 618-622.
\end{thebibliography}
\end{document}